\documentclass[3p,times]{elsarticle}

\usepackage[utf8x]{inputenc}
\usepackage{latexsym}

\usepackage{indentfirst}
\usepackage[english]{babel}
\usepackage{cancel}

\usepackage{lipsum}
\usepackage{array}
\usepackage{color} 
\usepackage{tabularx}
\usepackage{graphicx} 
\usepackage{amsmath}
\usepackage{amssymb}
\usepackage{amsfonts}	
\usepackage{moreverb}
\usepackage{dsfont}
\usepackage{tipa}
\usepackage{upgreek}
\usepackage{bm}
\usepackage{multirow}
\usepackage{soul}
\usepackage{ textcomp }
\usepackage[dvipsnames]{xcolor}
\usepackage{mathtools}
\usepackage{setspace}
\usepackage{algorithm}
\usepackage{algorithmic}
\usepackage{caption}
\usepackage{subcaption}
\usepackage{booktabs}
\allowdisplaybreaks

\newcommand\BibTeX{{\rmfamily B\kern-.05em \textsc{i\kern-.025em b}\kern-.08em
T\kern-.1667em\lower.7ex\hbox{E}\kern-.125emX}}

\usepackage{hyperref} 
\hypersetup{
    colorlinks=true,                          
    linkcolor=blue, 
    citecolor=red, 
    urlcolor=blue  } 

\newcommand{\ie}{{i.e.,}~}

\newcommand{\x}{\mbf{x}}

\newcommand{\mbf}[1]{\mathbf{#1}}			%

\renewcommand{\S}{\mathbf{S}}
\renewcommand{\u}{\mathbf{u}}
\newcommand{\w}{\mathbf{w}}

\newcommand{\q}{\mathbf{q}}
\newcommand{\F}{\mathbf{F}}
\newcommand{\f}{\mathbf{f}}

\renewcommand{\v}{\mathbf{v}}
\newcommand{\B}{\mathbf{B}}

\newcommand{\halb}{\frac{1}{2}}

\newcommand{\bdm}{\begin{displaymath}}
\newcommand{\edm}{\end{displaymath}}

\newcommand{\bea}{\begin{eqnarray} }
\newcommand{\eea}{\end{eqnarray} }

\newfont{\numerikEleven}{ecrm1000}
\newfont{\numerikTen}{cmss10}
\newfont{\numerikNine}{cmss9}
\newfont{\numerikEight}{cmss8}

\journal{Journal of Computational Physics}

\begin{document} 
\begin{frontmatter}
\title{High order numerical discretizations of the Einstein-Euler equations in the Generalized Harmonic formulation} 
\author[UniVR,UniTNmath]{Stefano Muzzolon}
\ead{stefano.muzzolon@univr.it}
\cortext[cor1]{Corresponding author} 

\author[UniTN]{Michael Dumbser}
\ead{michael.dumbser@unitn.it}

\author[UniTN]{Olindo Zanotti \corref{cor1}}
\ead{olindo.zanotti@unitn.it}

\author[UniVR]{Elena Gaburro}
\ead{elena.gaburro@univr.it}

\address[UniVR]{Department of Computer Science, University of Verona, Strada le Grazie 15, Verona, 37134, Italy}	
\address[UniTN]{Laboratory of Applied Mathematics, DICAM, University of Trento, via Mesiano 77, 38123 Trento, Italy} 
\address[UniTNmath]{Department of Mathematics, University of Trento, Via Sommarive 14, 38123 Trento, Italy}


\begin{abstract} \color[rgb]{0,0,0}
We propose two new alternative numerical schemes to solve the coupled Einstein-Euler equations in the Generalized Harmonic formulation.
The first one is a finite difference (FD) Central
Weighted Essentially Non-Oscillatory (CWENO) scheme on a traditional Cartesian mesh,
while the second one is an ADER (Arbitrary high order Derivatives)
discontinuous Galerkin (DG) scheme on 2D unstructured polygonal meshes. 
The latter, in particular, represents a 
preliminary step in view of a full 3D numerical relativity calculation on moving meshes.
Both schemes are equipped with a well-balancing (WB) property, which allows to preserve the equilibrium of \textit{a priori} known stationary solutions exactly at the discrete level. We validate our numerical approaches by successfully reproducing standard vacuum test cases, such as the robust stability, the linearized wave, and the gauge wave tests, as well as achieving long-term stable evolutions of stationary black holes, including Kerr black holes with extreme spin. Concerning the coupling with matter, modeled by the relativistic Euler equations, we perform some special relativistic Riemann problems, a classical test of spherical accretion onto a Schwarzschild black hole, as well as an evolution of  a perturbed non-rotating neutron star, demonstrating the capability of our schemes to operate also on the full Einstein-Euler system. Altogether, these results provide a solid foundation for addressing more complex and challenging simulations of astrophysical sources through DG schemes on unstructured 3D meshes.
\end{abstract}

\begin{keyword}
  Einstein field equations \sep relativistic Euler equations \sep Generalized Harmonic formulation \sep discontinuous Galerkin \sep non conservative \sep well-balancing 
 %
\end{keyword}
\end{frontmatter}


%
\section{Introduction} 
\label{sec.introduction}
When dealing with the numerical solution of the Einstein equations, both their mathematical formulation as a system of partial differential equations (PDEs) and their numerical discretization play a crucial role. At the moment, there are two
main successful  families of formulations of the Einstein equations that are routinely adopted for the simulation of astrophysical sources:  
those that are based 
on the 3+1 foliation of spacetime (see the monographies by~\cite{Alcubierre:2008},~\cite{Baumgarte2010},~\cite{Rezzolla_book:2013},~\cite{Gourgoulhon2012} for an overview), and those that belong to the Generalized Harmonic (GH) approach, which dates back to the origin of general relativity~\cite{DeDonder1927}. 
Within the 3+1 framework, after the seminal works by~\cite{Shibata95,Baumgarte99}, which gave birth to the so called BSSNOK formulation, 
there has been a steady effort over the years to obtain 
a better numerical control of the Einstein constraints. This has motivated the introduction of new versions such as Z4, CCZ4, Z4c, CFC~\cite{Bona:2003qn,Bona-and-Palenzuela-Luque-2005:numrel-book,Alic:2011a,Alic2013,Hilditch2013,Cordero-Carrion2008}. Within the GH framework, there has been continuous progress as well, which, not only allowed to perform
the  first stable numerical simulation of a binary black hole merger~\cite{Pretorius2005a}, but also made it possible to go beyond the rigidity of the original GH prescription regarding the gauge condition~\cite{Lindblom2006}.

If we look at the PDE structure of the two families of formulations, a fundamental feature emerges. While 3+1 formulations naturally arise as second-order (in space derivatives) PDE systems, GH formulations are easily written as first-order systems, which makes their hyperbolic nature much cleaner and more transparent to prove. This has strong implications on the choice of the numerical scheme that can be adopted for their solutions. In fact, one of the most promising numerical schemes, considered to be the frontier of numerical relativity (NR), is represented by discontinuous Galerkin (DG) schemes, which, notoriously, require a first-order formulation of the underlying PDEs. From this perspective, GH formulations are favored with respect to 3+1 second order formulations, 
and in fact DG schemes have already been successfully used for GH simulations, reaching a number of relevant scientific results~\cite{Kidder2017, Tichy2021, Tichy2023, Deppe2024}. In general, the implementation of high order numerical schemes for first-order formulations of the Einstein equations is a fertile field of research.
For example,  within the 3+1 framework, a new first-order version of the BSSNOK formulation that could evolve a black hole binary system with different high order numerical schemes, including DG, has been recently proposed in~\cite{DumbserZanottiPuppo2025, DumbserZanottiGaburroPeshkov2023}.

In this paper, we aim to give a contribution to the numerical discretization of the Einstein equations under the GH formulation, which consider various different perspectives.
First, we introduce a Central WENO Finite Difference (CWENO-FD) scheme on Cartesian grids, which has been successfully used within the first-order version of the BSSNOK formulation~\cite{DumbserZanottiPuppo2025}, but has never been adopted in the GH formulation.
Second, we introduce an ADER (Arbitrary high order Derivatives)
discontinuous Galerkin (DG) scheme on 2D unstructured polygonal meshes, which should be regarded as a preliminary implementation towards a truly 3D numerical relativity calculation under the GH framework.
In both cases we follow the choice of using a single numerical scheme for the whole PDE system, which includes the Einstein sector for the spacetime and the Euler sector for the matter. Besides guaranteeing a clean computational structure, this approach presents the crucial advantage of avoiding any sort of interpolation between different meshes, which may become necessary if different numerical schemes are used for the two sectors.
Moreover, both the CWENO and the DG schemes are enriched by a well-balancing property, which allows to preserve the equilibrium of \textit{a priori} known stationary (or quasi-stationary) configurations up to machine precision. This is an extra feature that can be activated according to user's needs and it follows closely the pragmatic approach proposed in~\cite{DumbserZanottiGaburroPeshkov2023}.

The plan of the paper is the following. In Sect.~\ref{sec.GHeqs} we recall the GH formulation of the Einstein equations; in Sect.~\ref{sec.numerical} we describe the two numerical schemes adopted in the paper, while Sect.~\ref{sec.tests} contains the validation of our approach through a wide set of numerical tests. Finally, Sect.~\ref{sec.conclusions} presents the conclusions of our work.

We work in a geometrized set of units, in which the speed of light and
the gravitational constant are set to unity, \ie $c=G=1$, and we assume $(-, +, +, +)$ as signature of the spacetime metric.
Latin indices from the first part of the
alphabet $a, b, c, \ldots $ denote spacetime indices ranging from 0 to 3, while Latin indices $i, j, k, \ldots$
are purely spatial, ranging from 1 to 3.
We also use the
Einstein summation convention of repeated indices.

\section{Basic features of the GH formulation}
\label{sec.GHeqs}

We assume  the  spacetime to be foliated through $\Sigma_t=const$ hypersurfaces as~\cite{Alcubierre:2008,Baumgarte2010,Rezzolla_book:2013}
\begin{equation}
	\label{eq:ds2_3p1}
	ds^{2} = g_{ab}\, dx^a\, dx^b=  -(\alpha^{2}-\beta_{i}\beta^{i}) dt^{2}+ 
	2 \beta_{i} dx^{i} dt + \gamma_{ij} dx^{i}dx^{j} \,,
\end{equation}
where $\alpha$ is the lapse function, $\beta^i$ is the shift vector and $\gamma_{ij}$ is the metric of the three dimensional space,
and where a natural observer is introduced with four velocity $n^a=(1/\alpha, -\beta^i/\alpha)$.
However, unlike 3+1 formulations where each vector or tensor is decomposed along directions parallel and perpendicular to $n^a$, in the Generalized Harmonic (GH) framework the splitting 
expressed by~\eqref{eq:ds2_3p1} does not represent a constitutive feature. 
In fact, the GH formulation of the Einstein equations involves the full four-dimensional metric tensor $g_{ab}$ and takes the first-order form~\cite{Deppe2024}

\begin{eqnarray}
	\label{eq:fosh gh metric evolution}
	\partial_t g_{ab}
	&=&\left(1+\gamma_1\right)\beta^k\partial_k g_{ab}
	-\alpha \Pi_{ab}-\gamma_1\beta^i\Phi_{iab}, \\
	\label{eq:fosh gh metric derivative evolution}
	\partial_t\Phi_{iab}
	&=&\beta^k\partial_k\Phi_{iab} - \alpha \partial_i\Pi_{ab}
	+ \alpha \gamma_2\partial_ig_{ab}
	+\frac{1}{2}\alpha n^c n^d\Phi_{icd}\Pi_{ab} \nonumber \\
	&+& \alpha \gamma^{jk}n^c\Phi_{ijc}\Phi_{kab}
	-\alpha \gamma_2\Phi_{iab},\\
	\label{eq:fosh gh metric conjugate evolution}
	\partial_t\Pi_{ab}
	&=&\beta^k\partial_k\Pi_{ab} - \alpha \gamma^{ki}\partial_k\Phi_{iab}
	+ \gamma_1\gamma_2\beta^k\partial_kg_{ab} \nonumber \\
	&+&2\alpha g^{cd}\left(\gamma^{ij}\Phi_{ica}\Phi_{jdb}
	- \Pi_{ca}\Pi_{db} - g^{ef}\Gamma_{ace}\Gamma_{bdf}\right) \nonumber \\
	&-&2\alpha \nabla_{(a}H_{b)}
	- \frac{1}{2}\alpha n^c n^d\Pi_{cd}\Pi_{ab}
	- \alpha n^c \Pi_{ci}\gamma^{ij}\Phi_{jab} \nonumber \\
	&+&\alpha \gamma_0\left(2\delta^c{}_{(a} n_{b)}
	- g_{ab}n^c\right)\mathcal{C}_c -\gamma_1\gamma_2 \beta^i\Phi_{iab}
	-16\pi \alpha \left(T_{ab} - \frac{1}{2}g_{ab}T\right),
\end{eqnarray}
where  we have used the definitions
\begin{eqnarray}
	\Phi_{iab}&=&\partial_i g_{ab}\,,\\
	\Pi_{ab} &=& -n^c\partial_cg_{ab}\,,
\end{eqnarray}
while $T_{ab}$ is the energy momentum tensor (see Eq.~\eqref{eq:Tmunu} below), $T$ is its trace, and $\delta_a^b$ is the Kronecker delta.
Regarding the notation, we recall that round brackets are used to indicate the symmetric part of a tensor, e.g. $A_{(ab)}=(A_{ab}+A_{ba})/2$, and that $\gamma^{ij}$ satisfies $\gamma^{ij}\gamma_{jk} = \delta^{i}_{k}$. In Eqs.~\eqref{eq:fosh gh metric evolution}--\eqref{eq:fosh gh metric conjugate evolution} three damping coefficients
are present:
$\gamma_0$, which damps the 1-index or gauge constraint
$\mathcal{C}_a=H_a+\Gamma_a$; $\gamma_1$, which controls the linear degeneracy of the
system; $\gamma_2$, which damps the 3-index constraint
$\mathcal{C}_{iab}=\partial_i g_{ab}-\Phi_{iab}$. Finally, $\Gamma_{abc}$ are the
spacetime Christoffel symbols of the first kind, \ie%
\begin{equation}
	\Gamma_{abc}=\frac{1}{2}(\partial_b g_{ac}+\partial_c g_{ab}-\partial_a g_{bc})\,,\qquad
\Gamma_a=g^{bc}\Gamma_{abc}\,,
\end{equation}
and $H_{a}$ is a \emph{gauge source function} which satisfies
\begin{equation}
	H_{a} = g_{ab} \nabla_{c} \nabla^{c} x^{b} = -\Gamma_{a}\,.
\end{equation}
During a typical simulation, we
will constantly monitor the so-called \emph{Einstein constraints}, given by 
\begin{eqnarray}
	\label{eqn.adm1}
	{\cal M}_a= \left( R_{ab}-\frac{1}{2}g_{ab}R - 8 \pi T_{ab} \right) n^b = 0\, ,
\end{eqnarray}
where
\begin{eqnarray}
	\label{eqn.ricci}
	R_{ab} = g^{cd} g^{ef} (\partial_{c} g_{ae}) (\partial_{d} g_{bf}) - g^{ec} g^{fd} \Gamma_{aef} \Gamma_{bcd} - \frac{1}{2} g^{cd} (\partial_{c} \partial_{d} g_{ab}) + \nabla_{(a}\Gamma_{b)} \, ,\qquad R = g^{ab}R_{ab}\,,
\end{eqnarray}
as well as the \emph{gauge constraints}, given by
\begin{equation}
	\label{eq:Cgauge}
	\mathcal{C}_a=H_a+\Gamma_a\,.
\end{equation}
It is also worth mentioning that the spatial components of $\Pi_{ab}$
are related to the usual extrinsic curvature tensor adopted within the classical 3+1 formulation through
\begin{equation}
	\label{eq:Pi-K}
	\Pi_{ij}= 2\,K_{ij} + \frac{\beta^{k}}{\alpha} \partial_{k}\gamma_{ij} - \frac{1}{\alpha}(D_{i}\, \beta_{j} + D_{j}\,\beta_{i}) \,,
\end{equation}
where $D_{i}$ represents the three-dimensional covariant derivative associated with the 3-metric $\gamma_{ij}$.
In our setup, the matter-field content is given by a perfect fluid only, hence with an 
energy-momentum tensor of the form
\begin{equation}
	\label{eq:Tmunu}
	T_{a b} = \rho h  u_a u_b + p g_{a b}\,,
\end{equation}
where $p$ is the gas pressure, $\rho$ is the rest mass density and $h$ is the specific (per unit mass) enthalpy, all of them measured in the comoving frame of the fluid with associated four velocity $u^{a}$. The latter is related to~$v^{i}$, the three velocity of the fluid as measured by the natural observer, by
\begin{equation}
    u^{i} = W \left(v^{i} - \frac{\beta^{i}}{\alpha}\right)\,,
\end{equation}
where $W$ is the Lorentz factor satisfying the relation
\begin{equation}
    W = \frac{1}{\sqrt{1-v^{2}}} = \alpha u^{0}\,,
\end{equation}
with $v^{2} = \gamma_{ij}v^{i}v^{j}$. The evolution equations for the fluid variables are the general relativistic Euler equations in conservation form~\cite{Banyuls97,DelZanna2007,Rezzolla_book:2013}
\begin{align}
\label{eqn.rho}
&\partial_t (\sqrt{\gamma}D)+\partial_i\left(\sqrt{\gamma}(\alpha v^i D - \beta^i D)\right)=0\,,\\ 
\label{eqn.S}
&\partial_t (\sqrt{\gamma} S_{j})+\partial_i\left(\sqrt{\gamma}(\alpha S^i_{\,\,j} - \beta^i S_j)\right)=\sqrt{\gamma}\left(\frac{\alpha}{2} S^{ik}\partial_{j}\gamma_{ik}  +
S_i  \partial_j \beta^{i} - E \partial_{j}\alpha  \right)\,,\\ 
\label{eqn.E}
&\partial_t (\sqrt{\gamma}E)+\partial_i\left(\sqrt{\gamma}(\alpha S^i - \beta^i E)\right)=\sqrt{\gamma}\left(\alpha S^{ij}K_{ij} - S^j \partial_{j}\alpha  \right)\,,
\end{align}
where $\gamma$ is the determinant of $\gamma_{ij}$, $S_{ij}$ is the spatial part of $T_{ab}$, and the remaining quantities in terms of the primitive physical variables read
\begin{equation}\label{eqn.cons.var}
    D = \rho W, \qquad S_{i} = \rho h W^{2} v_{i}, \qquad E = \rho h W^{2} - p\,.
\end{equation}
Using  a compact vector notation, the Einstein-Euler equations in GH formulation can be written as first order hyperbolic system of balance laws 
\begin{equation}
	\label{eqn.pde.mat.preview}
	\frac{\partial \u }{\partial t} + \nabla \cdot \F(\u)
	+ \B(\u) \cdot \nabla \u  
	= \S(\u),
\end{equation}
where the purely conservative term is only contributed by the matter part, with $\mathbf{F}(\u)=(\f^x(\u),\f^y(\u),\f^z(\u))$ the flux tensor, whereas the Einstein sector is entirely written in non-conservative form, with $\B(\u) = ( \B^x(\u), \B^y(\u), \B^z(\u) )$. 
The state vector $\u$ is composed of $54+5=59$ dynamical variables, \ie
10 for the metric tensor $g_{ab}$, 10 for the tensor $\Pi_{ab}$, 30 for the tensor $\Phi_{iab}$
, 4 for the vector $H_a$ and finally 
5 for the matter part. 
We remark that for the 54 variables of the Einstein sector the primitive and conservative variables coincide, while this is not the case for the 5 variables of the Euler sector.
Further details on the GH formulation can be found in~\cite{Alvi2002,Lindblom2006,Lindblom2009}.
%

\section{The numerical scheme}
\label{sec.numerical}
In our work, we consider two distinct families of numerical schemes: the Central WENO Finite Difference (CWENO-FD) scheme on Cartesian grids, and the ADER (Arbitrary high order Derivatives) discontinuous Galerkin (DG) scheme on 2D unstructured polygonal meshes. In what follows, we provide an overview of the main concepts characterizing each of the two methods, including the techniques used to obtain the well-balancing property, which allows to preserve an \textit{a priori} known equilibrium solution of the governing PDE system {\emph{exactly}} at the discrete level. 
We recall that, since their first introduction by~\cite{Bermudez1994}, well-balanced numerical methods have been extensively explored in the context of shallow water equations~\cite{leveque1998balancing, rebollo2003family, audusse2004fast, tang2004gas}, and, more generally, within the framework of hyperbolic balance laws~\cite{gosse2001well, perthame2001kinetic, BottaKlein, Gaburro2017CAFNonConf, Thomann2020, PareschiRey, castro2020well, berberich2021high,birke2024well}, eventually finding applications in numerical relativity as well~\cite{Gaburro2021WBGR1D, DumbserZanottiGaburroPeshkov2023, Balsara2024c, DumbserZanottiPuppo2025}. 
For a more theoretical foundations of these methods, one can refer for instance to~\cite{Castro2006, Pares2006, Castro2008}.

\subsection{A well-balanced path-conservative CWENO finite difference scheme}
\label{sec.WBCWENOFD}
Central WENO (CWENO) schemes were introduced in a series of seminal works for the numerical  solutions of hyperbolic systems of conservation laws in~\cite{Levy1999,Levy2000,Levy2001} and in the context of finite volume schemes in~\cite{cravero2016accuracy,SCR:CWENOquadtree,cravero2018cweno,tsoutsanis2021cweno}.
Their implementation is compatible with the more general framework of conservative WENO finite-difference methods~\cite{shuosher1,shuosher2,shu_efficient_weno,Balsara2024c}, for which some extensions to nonlinear hyperbolic systems involving non-conservative products were proposed in~\cite{Balsara2023,Balsara2024b,Balsara2024c}, allowing for a broader range of applications, including numerical relativity.
More recently, in the spirit of these studies, high order CWENO schemes have been successfully employed to tackle a novel first-order hyperbolic BSSNOK formulation, as presented in~\cite{DumbserZanottiPuppo2025}. The interested reader can find in there a detailed description of this approach, which we briefly summarize below.

\paragraph{Central WENO reconstruction}
We work with the discrete solution $\mathbf{u}_{i,j,k}$ of the system~\eqref{eqn.pde.mat.preview}, evaluated at the cell centers of a uniform Cartesian grid. For ease of explanation, we limit the description of  the procedure to the $x$ direction only. Hence, discrete points are placed at
\begin{equation}
	 x_i = x_L + \frac{1}{2} \Delta x + (i-1) \Delta x, 
	 \qquad 
	 \Delta x = \frac{x_R - x_L}{N_{x}}, 
\end{equation} 
where $x_L$ and $x_R$ are the left and right boundaries of the computational domain, while $N_x$ denotes the number of cells in that direction. 
We consider a large optimal stencil $\mathcal{S}_{opt} = \left\{ i{-}\frac{N}{2} , \, \dots, i{-}1, \, i, \, i{+}1, \dots, \, i{+}\frac{N}{2} \right\}$ , with $N > 2$ the degree of the reconstructed polynomial, and three smaller sub-stencils $\mathcal{S}_{L} = \left\{ i{-}2,\, i{-}1, \, i \right\},	\mathcal{S}_{C}=\left\{ i{-}1, \, i,\, i{+}1 \right\}, \mathcal{S}_{R} = \left\{ i, \, i{+}1, \, i{+}2 \right\}$ for polynomials of lower degree. 
Notice that, while we focus our description to the general case $N > 2$, the $N = 2$ case is treated analogously by employing only the side sub-stencils $\mathcal{S}_{L} = \left\{ i{-}1, \, i \right\}$ and $\mathcal{S}_{R} = \left\{ i, \, i{+}1 \right\}$.
We express a generic polynomial of degree $M = M(k)$ on $\mathcal{S}_k$ with respect to a set of polynomial basis functions  $\x \mapsto \psi_m(x)$ as
\begin{equation}
	P^M_k(x) = \sum \limits_{m=0}^M \psi_m(x) \, \hat u^m_k\,, \qquad \textrm{with} \qquad k \in \left\{ opt, L, C, R \right\},
\end{equation}
where the degrees of freedom $\hat u_k^m$ are obtained by imposing conservation over cells,
\begin{equation}
	\frac{1}{x_{j+\frac{1}{2}}-x_{j-\frac{1}{2}}} \int \limits_{x_{j-\frac{1}{2}}}^{x_{j+\frac{1}{2}}} P_k^M(x) dx = u_j, \qquad
	\forall j \in \mathcal{S}_k.  
\end{equation}
We compute an additional a priori unknown polynomial $P_0^N(x)$ from the optimal polynomial decomposition
\begin{equation}\label{eq.Popt}
P_{opt}^N(x) = \lambda_0 P_0^N(x) + \lambda_L P_L^2(x) + \lambda_C P_C^2(x) + \lambda_R P_R^2(x),
\end{equation}
where $\lambda_0, \lambda_L, \lambda_C, \lambda_R$ are arbitrary linear weights, which we set to $\lambda_0 = 10^8$, $\lambda_C = 10^4$, and $\lambda_L = \lambda_R = 1$, as in~\cite{dumbser2017central}.
We then use the standard WENO-type oscillation indicators $\sigma_{k}$ 
\begin{equation}
	\sigma_k = \sum_{\alpha=1}^{M} \int \limits_{x_{i-\halb} }^{x_{i+\halb }} \left( \frac{\partial^\alpha P^{M}_k(x)}{\partial x^\alpha} \right)^2 \Delta x^{2 \alpha-1} \, dx,\qquad \textrm{with}\qquad k \in \left\{ 0, L, C, R \right\},
\end{equation}
to carry out a non-linear WENO reconstruction in terms of  the polynomial $P_0^N(x)$ and of the lower degree polynomials, so that the final CWENO reconstruction is
\begin{equation}
	w_h(x) = \sum_k \omega_k P^{M}_k(x), \qquad \textnormal{ with } \qquad k \in \{0, L, C, R\},
\end{equation} 
where $\omega_{k}$ denote the nonlinear WENO weights 
\begin{equation}
	\label{WENOwr}
	{\omega}_k = \frac{\tilde{\omega}_k}{\sum_m \tilde{\omega}_m}, 
	\qquad 
	\textnormal{ with } 
	\qquad 
	\tilde{\omega}_k = \frac{\lambda_k}{\left( \sigma_k + \epsilon \right)^r}, 
\end{equation}
and where we set $r = 4$ and $\epsilon = 10^{-7}$.
See~\cite{DumbserZanottiPuppo2025} for further details.
To conclude, we evaluate $w_h(x)$ at the left and right interface of each interval defining $w^{\mp}_{i \pm \frac{1}{2},j,k} = w_h \left( x_{i\pm\frac{1}{2}}\right)$, and we use these reconstructed states to produce a high order accurate approximation of the first spatial derivative
\begin{equation}
	\partial_x \u_{i,j,k} = \frac{	\w^-_{i + \frac{1}{2},j,k} -	\w^+_{i - \frac{1}{2},j,k}}{\Delta x}.
\end{equation} 
This type of reconstruction is repeated for the discrete derivatives along the $y$ and the $z$ directions, to obtain the discrete spatial gradient $\nabla \u_{i,j,k}$, and it is also applied to the point values of the fluxes, again in each direction. 
Concerning the variables to be reconstructed, we mention that in principle  one could try to adopt characteristic variables using the available form of the analytic eigenvectors.	
However, due to the excessive computational cost, it is virtually impossible 
to project repeatedly the entire system into that eigenspace. We therefore simply employ componentwise reconstruction in conservative variables. 
\paragraph{CWENO-FD discretization}
The semi-discrete path-conservative FD scheme for system ~\eqref{eqn.pde.mat.preview} then reads
\begin{eqnarray}
	\label{semi-discrete-scheme}
	\frac{d \u_{i,j,k} }{d t} &=& 
	- \frac{\f^x_{i+\halb,j,k}-\f^x_{i-\halb,j,k}}{\Delta x}   
	- \frac{\f^y_{i,j+\halb,k}-\f^y_{i,j-\halb,k}}{\Delta y}   
	- \frac{\f^z_{i,j,k+\halb}-\f^z_{i,j,k-\halb}}{\Delta z}  
	- \B(\u_{i,j,k}) \cdot \nabla \u_{i,j,k}    \nonumber \\
	&& 
	- \frac{\mathbf{D}^x_{i+\halb,j,k}+\mathbf{D}^x_{i-\halb,j,k}}{\Delta x}  -
	\frac{\mathbf{D}^y_{i,j+\halb,k}+\mathbf{D}^y_{i,j-\halb,k}}{\Delta y}  -
	\frac{\mathbf{D}^z_{i,j,k+\halb}+\mathbf{D}^z_{i,j,k-\halb}}{\Delta z}   + \S(\u_{i,j,k} ),
\end{eqnarray}
where $\partial_{x} \mathbf{u}_{i,j,k}$ is determined using the CWENO reconstruction reported above, and the numerical fluxes are obtained through a Rusanov-type Riemann solver. The jump terms involving $\mathbf{D}^x$, which do not appear in traditional central FD schemes, affect the order of accuracy and the robustness of the overall method, see~\cite{Balsara2024c}.
For the Rusanov-type scheme that we adopt, they are defined as
\begin{equation}\label{eq.Dterms}
\mathbf{D}^x_{i + \frac{1}{2},j,k} =  \frac{1}{2}  \bar{\mathbf{B}}^x_{i+\frac{1}{2},j,k} \cdot   
\left( \mathbf{w}^{+}_{i + \frac{1}{2},j,k} - \mathbf{w}^{-}_{i + \frac{1}{2},j,k} \right)  \,,\\ 
\end{equation}
where we use a simple midpoint rule for the $\bar{\mathbf{B}}^x$ term, namely  $\bar{\mathbf{B}}^x_{i+\frac{1}{2},j,k}=\mathbf{B}^x\left(\frac{1}{2}(\mathbf{w}^{+}_{i + \frac{1}{2},j,k}+\mathbf{w}^{-}_{i + \frac{1}{2},j,k})\right)$. The numerical flux reads
\begin{equation}
	\label{Rusanov}
	\f^x_{i+\halb,j,k}=\frac{1}{2}\left(\f^{x,-}_{i + \halb,j,k} + \f^{x,+}_{i + \halb,j,k}   \right)-\frac{1}{2}\lambda^{\max}_{i + \halb,j,k}\left( \w^{+}_{i + \halb,j,k} - \w^{-}_{i + \halb,j,k} \right), 
\end{equation}
with $\lambda^{\max}_{i + \halb,j,k}$ an estimate of the maximum wave speed at the interface.
Analogous formulae are used in the $y$ and $z$ directions, respectively.

For the temporal integration, we employ classical or TVD Runge-Kutta schemes of sufficiently high order of accuracy, consistent with the spatial discretization.

\paragraph{Well-balancing of the final CWENO-FD scheme}
To make the numerical scheme well-balanced, we follow the approach described in~\cite{DumbserZanottiGaburroPeshkov2023}, which is in turn inspired by~\cite{Ghosh2016,PareschiRey,berberich2021high}.
The procedure is in principle quite simple: we subtract the equilibrium solution from the governing PDE system, thus removing the discretization errors introduced by the numerical scheme in the vicinity of the equilibrium.
In other words, if $\mathbf{u}_{e} = \mathbf{u}_{e}(x,y,z)$ denotes any \textit{a priori} known stationary equilibrium solution of the PDE~\eqref{eqn.pde.mat.preview}, so that $\partial_{t} \mathbf{u}_{e} = 0$, we subtract Eq. ~\eqref{eqn.pde.mat.preview} evaluated at $\mathbf{u}_{e}$ from its general form, thus obtaining the \textit{augmented system}
\begin{eqnarray}
	\frac{\partial \mathbf{u} }{\partial t} +\nabla \cdot \mathbf{F}(\mathbf{u}) 
	- 
	\nabla \cdot \mathbf{F}(\mathbf{u}_{e}) 
	+
	\mathbf{B}(\mathbf{u}) \cdot \nabla \mathbf{u}
	-
	\mathbf{B}(\mathbf{u}_e) \cdot \nabla \mathbf{u}_{e} 
	&=& \mathbf{S}(\mathbf{u}) - \mathbf{S}(\mathbf{u}_e), \\
	\frac{\partial \u_e }{\partial t}  & = & 0. 
\end{eqnarray}
By introducing a new vector variable $\tilde{\u} = [\u, \u_{e}]^{T}$, which is evolved by the code, and suitable extensions of the operators involved, the above system can be cast into the compact form of the original PDE, namely
\begin{equation}\label{eqn.pde.mat.preview.aug}
	\frac{\partial \tilde{\u}}{\partial t} + \nabla \cdot \tilde{\F}(\tilde{\u}) + \tilde{\B}(\tilde{\u}) \cdot \nabla \tilde{\u}   
	= \tilde{\S}(\tilde{\u}),
\end{equation} 
with the augmented flux tensor, non-conservative product, and source term given by
\begin{equation}
	\tilde{\F} = 
	\left( \begin{array}{c} 
		\mathbf{F}(\mathbf{u}) - \mathbf{F}(\mathbf{u}_{e}) \\ 
		0 
	\end{array} \right), \qquad 
	\tilde{\B}(\tilde{\u}) \cdot \nabla \tilde{\u} = 
	\left( \begin{array}{c} 
		\B(\u) \cdot \nabla \u  - \B(\u_e) \cdot \nabla \u_e \\ 
		0  
	\end{array} \right), \qquad 
	\tilde{\S} = 
	\left( \begin{array}{c} 
		\mathbf{S}(\mathbf{u}) - \mathbf{S}(\mathbf{u}_{e}) \\ 
		0 
	\end{array} \right).                  
\end{equation} 
In order to obtain a well-balanced numerical dissipation, we emphasize that it is necessary to modify the numerical flux~\eqref{Rusanov} for the system~\eqref{eqn.pde.mat.preview.aug} as follows
\begin{equation}
	\label{Rusanov.wb}
	\tilde{\f}^x_{i+\halb,j,k}=\frac{1}{2}\left(\tilde{\f}^{x,-}_{i + \halb,j,k} + \tilde{\f}^{x,+}_{i + \halb,j,k}   \right)-\frac{1}{2}\lambda^{\max}_{i + \halb,j,k} \tilde{\mathbf{I}} \left( \tilde{\w}^{+}_{i + \halb,j,k} - \tilde{\w}^{-}_{i + \halb,j,k} \right), 
\end{equation}
with 
\begin{equation}  
	\tilde{\mathbf{I}} = \left(  
	\begin{array}{c c c}  
		\mathbf{I} & \vline & -\mathbf{I} \\  
		\hline   
		0 &\vline & 0  
	\end{array}  
	\right)
	\label{eqn.wb.id}
\end{equation}
being the so-called \textit{well-balanced identity matrix}, while $\mathbf{I}$ is the classical identity matrix. 
It is then obvious that for solutions that satisfy $\tilde{\u} = [\u_e,\u_e]$ we get by construction 
\begin{equation}
	\frac{d \tilde{\u}_{i,j,k} }{d t} = 0 
\end{equation}  
exactly at the semi-discrete level. For further details, see ~\cite{DumbserZanottiGaburroPeshkov2023}.

\subsection{A well-balanced path-conservative ADER discontinuous Galerkin scheme}
\label{sec.WBADERDG}
The ADER approach was first introduced in~\cite{toro2001towards} as a fully-discrete one-step procedure to achieve arbitrary high order of accuracy in both space and time. As shown  in~\cite{QiuDumbserShu, dumbser_jsc, GassnerDumbserMunz}, it can be framed within the context of discontinuous Galerkin (DG) finite element methods, leading to the so called ADER-DG schemes. Its modern formulation~\cite{DumbserEnauxToro,dumbser2008unified}, which consists in using a local predictor-corrector technique, has demonstrated remarkable effectiveness and versatility on many different PDE systems~\cite{DumbserZanotti, Dumbser2009a, Dumbser2010, ADERNSE, BalsaraMultiDRS, BalsaraDivB2015, ADERGRMHD, Zanotti2015d}, including some applications to general relativity~\cite{Dumbser2017strongly, Dumbser2020GLM, Gaburro2021WBGR1D,DumbserZanottiPeshkov2024}. In addition, it has also been extended to the Arbitrary-Lagrangian-Eulerian (ALE) framework on general unstructured polygonal meshes~\cite{Lagrange2D,  Lagrange3D, ALELTS2D, ALEDG, Gaburro2020, Gaburro2020Arepo}. The latter is the setting that we assume here, although in the present paper we limit our attention  to a fixed mesh configuration, leaving the much more challenging case of moving meshes with topology changes~\cite{Springel,Gaburro2020Arepo,gaburro2021high} to a future work.

\begin{figure}[!htbp]
\begin{center}
\includegraphics[width=0.8\textwidth, keepaspectratio]{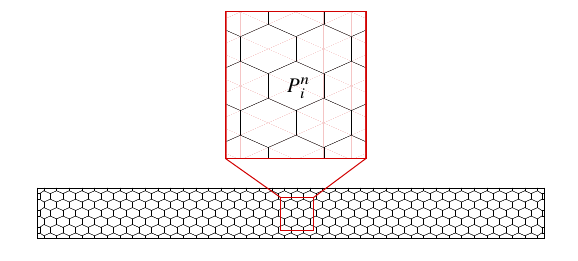}
\caption{
An example of a rectangular 2D unstructured mesh covered by non-overlapping polygons $P_{i}^{n}$. Above, a zoomed-in portion highlights the polygonal structure.
}
\label{fig.unstructured.mesh.example}
\end{center}	
\end{figure}

\paragraph{Data representation}
We consider a 2D domain covered by an unstructured mesh of non-overlapping polygons $P_i^{n}$, $i = 1, \dots, N_{P}$ (see Figure~\ref{fig.unstructured.mesh.example} for a representative illustration), that could be in principle rearranged at each time step $t^{n}$, where the vector variable $\mathbf{u}$ for Eq.~\eqref{eqn.pde.mat.preview} is represented by a piecewise polynomial of degree $N$, 
\begin{equation}\label{eqn.uh}
\mathbf{u}_h^n(\x,t^n) = \sum \limits_{\ell=0}^{\mathcal{N}-1} \varphi_\ell(\x,t^n) \, \hat{\mathbf{u}}^{n}_{\ell,i}, \qquad \mathbf{x} \in P_i^n,
\end{equation}
with $\x \mapsto\varphi_{\ell}(\mathbf{x}, t^n)$ local modal spatial basis functions, and $\mathcal{N}= \mathcal{L}(N,2)$ where $\mathcal{L}(N,d) = \frac{1}{d!} \prod \limits_{m=1}^{d} (N+m)$.
In practice, to determine the total number of elements $N_{P}$ in the typical case of a rectangular domain, we fix a desired approximate number of elements in the $x$ direction $N_{x}$, and we then compute an optimal number of elements in the $y$ direction $N_{y}$, to ensure the resulting mesh is as uniform as possible. Consequently, in the numerical tests, we will use $N_{x}$ as an indicative measure of the mesh resolution as it directly controls the refinement process. 

\paragraph{The space-time predictor}
\label{sec.predictor}
The main ingredient of the method is a {\emph{local}} procedure consisting in solving a weak form of the governing PDE system to obtain a {\emph{predictor}} $\mathbf{q}_h^n$ of the solution for each element, which is a space-time polynomial of degree $N$ described by
\begin{equation}
\mathbf{q}_h^n(\mathbf{x}, t) = \sum_{\ell=0}^{\mathcal{Q}-1} \theta_{\ell}(\mathbf{x}, t) \hat{\mathbf{q}}_{\ell}^n, \qquad (\mathbf{x}, t) \in C_i^n,
\end{equation}
where $(\x,t) \mapsto \theta_{\ell}(\mathbf{x}, t)$ are modal space-time basis functions, $\mathcal{Q}= \mathcal{L}(N,3)$, and ${C_i^n}$ is the space-time control volume bounded by the elements $P_i^{n}$ and $P_i^{n+1}$.
In more detail, we multiply the PDE system~\eqref{eqn.pde.mat.preview} by time-dependent test functions $(\x,t) \mapsto \theta_{k}(\mathbf{x}, t)$ and we integrate over $ C_i^n$, obtaining
\begin{equation} \begin{aligned} 
\label{eqn.predictor.1}
\int_{C_i^n} \theta_k(\mathbf{x} ,t) \frac{\partial \mathbf{q}_h^n}{\partial t} \, d \mathbf{x} \, d t \,
& +\int_{C_i^n} \theta_k(\mathbf{x} ,t) \left( \nabla \cdot \mathbf{F}(\mathbf{q}_h^n) + \mathbf{B}(\mathbf{q}_h^n) \cdot \nabla \mathbf{q}_h^n  \right) \, d\mathbf{x} \, d t \\
& = \int_{C_i^n} \theta_k(\mathbf{x} ,t) \mathbf{S}(\mathbf{q} _h^n) \, d \mathbf{x}  \, d t\,, 
\quad \forall k \in \left\{0, \dots, \mathcal{Q}-1\right\}.
\end{aligned} \end{equation}
Since in a space-time DG approach in general $\q_h(\x,t^{n}) \neq \u_h(\x,t^n)$, we need to account for the jump at time $t^n$ in the sense of distributions. Making use of the causality principle, upwinding in time leads to the following weak form in space-time

{
\begin{equation} \begin{aligned}
			\label{eqn.predictor.2}
			& \int_{P_i^{n}}   \theta_k(\mathbf{x},t^{n}) \left(  \mathbf{q}_h^n(\mathbf{x},t^{n})  -  \u_h^n(\mathbf{x},t^{n}) \right)  \,d \mathbf{x} \, 
			+ \int_{C_i^{n}} \theta_k(\mathbf{x},t) \frac{\partial }{\partial t}  \mathbf{q}_h^n(\mathbf{x},t) \,d \mathbf{x} \,d t
			 + \ \int_{C_i^{n}}   \theta_k(\mathbf{x},t) \nabla \cdot \mathbf{F}(\mathbf{q}_h^n)  \,d \mathbf{x} \,d t = 
			 \\ &   
			\int_{C_i^{n}}  \theta_k(\mathbf{x},t)
			\left( \mathbf{S}(\mathbf{q}_h^n) - \mathbf{B}(\mathbf{q}_h^n) \cdot \nabla \mathbf{q}_h^n \right) \,d \mathbf{x} \,d t\,, 
			\quad \forall k \in \left\{0, \dots, \mathcal{Q}-1\right\},
\end{aligned} 
\end{equation}
}%
which translates to an element-local nonlinear system for the degrees of freedom $\hat{\mathbf{q}}_{\ell}^n$, that we solve via a discrete Picard iteration~\cite{DumbserEnauxToro,Zanotti2015}, whose convergence was proven in~\cite{BCDGP20}.

\paragraph{The corrector step}
We start by rewriting Eq.~\eqref{eqn.pde.mat.preview} in a space-time divergence form,  
\begin{equation}\label{eqn.pde.divform}
\begin{aligned}
\tilde \nabla \cdot \tilde{\mathbf{F}}(\mathbf{u}) + \tilde{\mathbf{B}}(\mathbf{u}) \cdot \tilde{\nabla} \mathbf{u} = \mathbf{S}(\mathbf{u}), \quad \text{ with }    
\tilde{\mathbf{F}} = (\mathbf{u}, \mathbf{F}), \, \tilde{\mathbf{B}} = (\mathbf{0}, \mathbf{B}),  \text{ and }  
\tilde \nabla  = \left( \partial_t, \, \partial_\mathbf{x} \right)^T.
\end{aligned}
\end{equation}
We then multiply it by so-called moving test functions $(\x,t) \mapsto \tilde{\varphi}_{k}(\x,t)$, 
which coincide with the basis function of Eq.~\eqref{eqn.uh} at $t=t^n$ and at $t=t^{n+1}$, 
i.e. $\tilde{\varphi}_k(\x,t^n) = \varphi_k(\x,t^n)$ and $\tilde{\varphi}_k(\x,t^{n+1}) = \varphi_k(\x,t^{n+1})$ and automatically  adapt to an eventual moving mesh framework. 
Next, we integrate over each $C_i^n$ applying the Gauss theorem to the flux-divergence term and splitting the non-conservative products into volume and surface contribution, to obtain 
\begin{equation}\label{eqn.corrector}
\begin{aligned}
&\int_{P_i^{n+1}} \tilde{\varphi}_k  \mathbf{u}_h(\mathbf{x},t^{n+1}) \, d\mathbf{x}   
=  \int_{P_{i}^n} \tilde{\varphi}_k  \mathbf{u}_h(\mathbf{x},t^{n})  \, d\mathbf{x} 
- \sum_{j=1}^{\lvert \mathcal{V}_i^{n} \rvert} \int_{\partial C_{ij}^n} \tilde{\varphi}_k  \mathcal{D}(\mathbf{q}_h^{n,-},\mathbf{q}_h^{n,+}) \cdot \mathbf{\tilde{n}} \, dS \\
& + \int_{C_i^{n}} \tilde \nabla \tilde{\varphi}_k \cdot \tilde{\mathbf{F}} (\mathbf{q}_h)\, d\mathbf{x} dt  
+ \int_{C_i^{n}}  \tilde{\varphi}_k(\mathbf{x},t) \left( \mathbf{S}(\mathbf{q}_h^n) - \mathbf{B}(\mathbf{q}_h^n) \cdot \nabla \mathbf{q}_h^n \right) \, d\mathbf{x} \, d t\,, 
\quad \forall k \in \left\{0, \dots, \mathcal{N}-1\right\},
\end{aligned}
\end{equation}
where $\mathcal{V}_i^{n}$ is the set of neighbors of $P_i^{n}$, $\mathbf{\tilde{n}} = (\tilde{n}_{t}, \, \mathbf{\tilde{n}_{x}})$ denotes the outward pointing spacetime unit normal vector on the faces composing the boundary of $C_{i}^{n}$ , the $\tilde{\mathbf{F}}$-term is simply the physical space-time flux evaluated at the discrete solution $\mathbf{q}_h$ and the $\mathcal{D}$-term is a path-conservative Rusanov-type fluctuation that reads 
\begin{equation}
	\mathcal{D}(\mathbf{q}_h^{n,-},\mathbf{q}_h^{n,+}) \cdot \mathbf{\tilde{n}} = 
	\frac{1}{2} \left( 
	\tilde{\mathbf{F}} (\mathbf{q}_h^+) +  	\tilde{\mathbf{F}} (\mathbf{q}_h^-) \right) \cdot \mathbf{\tilde{n}} 
	+ \frac{1}{2} \bar{\mathbf{B}} \cdot \left( \mathbf{q}_h^+ - \mathbf{q}_h^- \right) 
	- \frac{1}{2} s_{\max}  \left( \mathbf{q}_h^+ - \mathbf{q}_h^- \right).
	\label{eqn.pc.rusanov} 
\end{equation}
Here, $\q_h^{n,-}$ and $\q_h^{n,+}$ are the boundary-extrapolated data from the predictor in each shared lateral surface $\partial C_{ij}^{n}$, $ s_{\max}$ is an estimate for the maximum signal speed at the cell interface and the average matrix $\bar{\mathbf{B}}$ is
computed according to a path integral along a straight-line segment path as follows: 
\begin{equation}
	\bar{\mathbf{B}} = \int \limits_0^1 \tilde{\mathbf{B}}(\psi(\xi)) \cdot \mathbf{\tilde{n}} \, d\xi, 
	\qquad 
	\psi = \mathbf{q}_h^- + \xi \left( \mathbf{q}_h^+ - \mathbf{q}_h^- \right).
\end{equation}
Eq.~\eqref{eqn.corrector} provides the final update from $t^n$ to $t^{n+1}$, since the unknown $\mathbf{u}_h(\mathbf{x},t^{n+1})$ can be computed from  the solution at the previous time step, $\mathbf{u}_h(\mathbf{x},t^{n})$ by exploiting the predictor $\mathbf{q}_h(\mathbf{x},t)$ when computing the various integrals appearing in the formulation. 

\paragraph{Well-balancing of the final ADER-DG scheme}
For the well-balancing of the scheme, we refer to the approach recently precisely described in~\cite{GaburroALEWB} and references therein. The underlying idea is to consider three groups of quantities: the equilibrium solution $\mathbf{u}_{e}$, \textit{a priori} known in space and time, which exactly satisfies the PDE system~\eqref{eqn.pde.mat.preview}; the {\emph{fluctuation}} $\mathbf{u}_{f}$, which accounts for the deviations from the equilibrium profile and is also represented following Eq.~\eqref{eqn.uh}; and the complete numerical solution $\mathbf{u}_{h}$, which corresponds to the variable involved in the original (non well-balanced) version of the scheme, that however here takes the form
\begin{equation}\label{eq.recover.fluct}
	\mathbf{u}_{h}(\mathbf{x}, t) = \mathbf{u}_{e}(\mathbf{x}, t) + \mathbf{u}_{f}(\mathbf{x}, t),
\end{equation}
accounting for both the contributions of the equilibrium and its fluctuations. Based on these definitions, we proceed as follows. 

We first subtract the predictor equation~\eqref{eqn.predictor.2} evaluated at the equilibrium state $\u_e$, to Eq.~\eqref{eqn.predictor.2} itself, obtaining a system for the unknown fluctuation 
\begin{equation}\label{eq.recover.predict.fluct}
	\q_{f}(\x, t) = \q_{h}(\x, t) - \u_{e}(\x, t),
\end{equation}
namely
\begin{equation} \begin{aligned}
		\label{eqn.predictor.fluct}
		& \int_{P_i^{n}}   \theta_k(\x,t^{n}) \,  \q_f^n(\x,t^{n})  \,d \x \,
		+ \int_{C_i^{n}} \theta_k(\x,t) \frac{\partial \q_f^n}{\partial t}  \,d \x \,d t = \int_{P_i^{n}}   \theta_k(\x,t^{n})  \u_f^n(\x,t^{n})  \,d \x \, \\
		&- \int_{C_i^{n}}   \theta_k(\x,t) \left[ \nabla \cdot \F(\q_h^n) - \nabla \cdot \F(\u_e) \right] \,d \x \,d t   \\
		&+\int_{C_i^{n}}  \theta_k(\x,t) \left[
		\left( \S(\q_h^n) - \B(\q_h^n)  \cdot \nabla \q_h^n \right) - \Bigl( \S(\u_e) - \B(\u_e) \cdot \nabla \u_e \Bigr) \right] \,d \x \,d t \,, 
		\quad \forall k \in \left\{0, \dots, \mathcal{Q}-1\right\},
	\end{aligned} 
\end{equation}
which we solve using a fixed point iteration as before. 	
Notice that the systematic evaluation of differences between terms associated with $\q_h$ and $\u_e$ significantly reduces the numerical errors when $\u_h$ is a small perturbation of $\u_e$, which indeed represents the natural case of application of well-balancing techniques. Moreover, if $\u_f = \mathbf{0}$ then $\q_f = \mathbf{0}$, meaning that the resulting predictor is \textit{exact} on equilibria.
	
To complete our well-balanced scheme, we then use the same subtraction strategy for the corrector equation~\eqref{eqn.corrector} to derive the following update formula for the degrees of freedom of the fluctuation $\u_f$, again based on the relation~\eqref{eq.recover.fluct},
\begin{equation}\label{eqn.corrector.fluct}
	\begin{aligned}
		&\int_{P_i^{n+1}} \tilde{\varphi}_k  \, \u_f(\x,t^{n+1})  \, d\mathbf{x}   
		=  \int_{P_{i}^n} \tilde{\varphi}_k \, \u_f(\x,t^{n}) \, d\mathbf{x}  + \int_{C_i^{n}} \tilde \nabla \tilde{\varphi}_k \cdot \left( \tilde{\mathbf{F}} (\mathbf{q}_h) - \tilde{\mathbf{F}} (\u_e) \right) \, d\mathbf{x} dt \\
		&- \sum_{j=1}^{\lvert \mathcal{V}_i^{n} \rvert} \int_{\partial C_{ij}^n} \tilde{\varphi}_k \left(  \mathcal{D}(\mathbf{q}_h^{n,-},\mathbf{q}_h^{n,+}) \cdot \mathbf{\tilde{n}} - \mathcal{D}(\u_e^{-},\u_e^{+}) \cdot \mathbf{\tilde{n}} \right) \, dS\\
		&+ \int_{C_i^{n}}  \tilde{\varphi}_k(\mathbf{x},t) \left[ \left( \mathbf{S}(\mathbf{q}_h^n) - \mathbf{B}(\mathbf{q}_h^n) \cdot \nabla \mathbf{q}_h^n \right)- \Bigl( \mathbf{S}(\u_e) - \mathbf{B}(\u_e) \cdot \nabla \u_e \Bigr) \right] \, d\mathbf{x} \, d t\,,
		\quad \forall k \in \left\{0, \dots, \mathcal{N}-1\right\},
	\end{aligned}
\end{equation}
where, for smooth equilibria, the boundary extrapolated values of $\u_e$ coincide. It is important to underline that the predictor $\q_h$ to be used inside Eq.~\eqref{eqn.corrector.fluct} should be obtained by summing up the equilibrium solution to the computed predictor of the fluctuation, specifically by rearranging Eq.~\eqref{eq.recover.predict.fluct}.
	
As a final remark, we emphasize that the well-balancing properties described above for the predictor step are inherited by the corrector update, and consequently, the overall scheme is exact on equilibria and minimizes the numerical errors in presence of small fluctuations.
\paragraph{A posteriori subcell finite volume limiter}
We need to mention that high order discontinuous Galerkin (DG) schemes are \textit{linear} in the sense of Godunov~\cite{godunov} and, consequently, are affected by the Gibbs phenomenon, producing detrimental oscillations in presence of discontinuities. 
A variety of techniques have been proposed over the years to address this issue, where a representative though non-exhaustive selection includes slope~\cite{cockburn1998runge}, moment~\cite{krivodonova2004shock} and damping~\cite{lu2021oscillation} limiters, artificial viscosity~\cite{persson2006sub} and WENO-type~\cite{qiu2004hermite, qiu2005runge} approaches, spectral filtering~\cite{Radice2011}, gradient-based procedures~\cite{kuzmin2020gradient}, and \textit{a posteriori} subcell finite volume (FV) limiters~\cite{DGLimiter1, DGLimiter2, Zanotti2015d, DGLimiter3, Gaburro2020Arepo, gaburro2021posteriori, gaburro2024discontinuous, GaburroALEWB}. We adopt the last approach, which is inspired by the Multi-dimensional Optimal Order Detection (MOOD) method~\cite{clain2011high, loubere2014new} and has already been effectively applied in the context of numerical relativity, for instance in~\cite{Deppe2022, Dumbser2017strongly, ADERGRMHD}. In the following, we provide only a summary of the key aspects related to our scheme, while referring to the cited literature for more comprehensive and detailed descriptions.

We first apply the ADER-DG scheme outlined above over the entire computational domain in order to obtain a \textit{candidate} solution at time $t^{n+1}$, denoted by $\u^{n+1,*}_{h/f}$, where the subscript refers either to the standard discrete solution or to the fluctuation when the well-balancing feature is activated.

For each polygon $P_{i}^{n}$ we introduce a subtriangulation $\mathcal{T}(P_{i}^{n})$ composed of $\lvert \mathcal{V}_i^{n} \rvert$ triangles, each of which is further divided into $N^{2}$ smaller subtriangles $s_{i, \, \alpha}^{n}$, with $\alpha = 1, \dots, \lvert \mathcal{V}_i^{n} \rvert \cdot N^{2}$. We then define the quantities $\v_{i, \, \alpha}^{n+1}$ as the projections of $\u^{n+1,*}_{h}$ onto $s_{i, \, \alpha}^{n+1}$ via
\begin{equation}\label{eq.proj.subtr}
	\v_{i, \, \alpha}^{n+1} = \frac{1}{\lvert s_{i, \, \alpha}^{n+1} \rvert} \int_{s_{i, \, \alpha}^{n+1}} \u^{n+1,*}_{h}(\x, t^{n+1})\, d\x \,, \quad \forall \alpha \,.
\end{equation}%
At this stage, we check both the values of $\u^{n+1,*}_{h}$ on each vertex of $s_{i, \, \alpha}^{n+1}$ and $\v_{i, \, \alpha}^{n+1}$ itself against several admissibility criteria of different nature: numerical ones, ensuring the absence of floating point errors; physical ones, requiring positive densities and pressures, and that no superluminal velocities arise; and bounds-preserving ones, enforcing a relaxed discrete maximum principle (DMP).
Notice that in the well-balanced case we anyway check the solution $\u^{n+1,*}_{h}$, obtained via Eq.~\eqref{eq.recover.fluct}, and not the fluctuation  $\u^{n+1,*}_{f}$.
On any polygon where the candidate solution meets the required conditions, we directly accept it as is by setting $\u^{n+1}_{h/f} = \u^{n+1,*}_{h/f}$. On the other hand, if this is not the case, we mark the candidate solution as \textit{troubled}, we discard it, and recompute it locally using a more robust approach, specifically a second order ADER-ENO finite volume (FV) scheme, applied this time to the subtriangulation. This detail is essential to retain as much as possible the accuracy provided by the DG method.
More precisely, we use a projection analogous to Eq.~\eqref{eq.proj.subtr} to compute the cell averages $\v_{i, \alpha}^{n}$ from the available solution at time $t^{n}$, and we evolve them with the finite volume method to obtain the new $\v_{i, \alpha}^{n+1}$. This allows us to reconstruct a DG polynomial through
\begin{equation}\label{eq.reconstr.fv.dg}
	\int_{s_{i, \, \alpha}^{n}} \u^{n+1}_{h}(\x, t^{n+1})\, d\x = \int_{s_{i, \, \alpha}^{n}} \v^{n+1}_{h}\, d\x \,, \quad \forall \alpha \,,
\end{equation}
while imposing conservation over the whole $P_{i}^{n+1}$ by means of
\begin{equation}\label{eq.conserv.fv.dg}
	\int_{P_{i}^{n+1}} \u^{n+1}_{h}(\x, t^{n+1})\, d\x = \int_{P_{i}^{n+1}} \v^{n+1}_{h}\, d\x \,,
\end{equation}
thus closing the algorithm.

We remark that these final steps do not destroy the well-balancing property of the scheme, provided that the cell averages $\v_{i, \, \alpha}^{n}$ are divided into equilibrium and fluctuation parts, and that a well-balanced finite volume scheme is used to evolve the fluctuations $\v_{f, \, i, \, \alpha}^{n}$.

\section{Numerical tests}
\label{sec.tests}
In this Section, we present an extensive set of standard test cases in numerical relativity to assess the accuracy and robustness of the two numerical methods proposed in Sect.~\ref{sec.numerical}, namely the CWENO Finite Difference and ADER discontinuous Galerkin schemes, applied to the first-order hyperbolic Generalized Harmonic (GH) formulation of the Einstein-Euler system.
In particular, the vacuum spacetime setups are taken from the widely used and well-established benchmark tests provided in~\cite{Alcubierre2004}.
In what follows, we briefly describe each test case and analyze the effectiveness of our numerical approaches, especially in achieving long-term stable simulations, which, unless otherwise specified,
are all performed up to the canonical final time of $t =1000$. 
To this end, we monitor the time evolution of the Einstein and gauge constraints given by Eqs.~\eqref{eqn.adm1} and~\eqref{eq:Cgauge}, and referring to their components as $M0,M1,M2,M3$, and $C0,C1,C2,C3$, respectively. We remark that, throughout all the tests, we show only the distinct components of the constraints when some coincide, and we omit from the plots any components that vanish identically.
We also mention that we consider the components of the gauge source function $H_{a}$ to be fixed once specified on the initial data, even though one could provide evolution equations also for them if necessary~\cite{Pretorius2005a, Pretorius2005b, lindblom2008gauge}. Moreover, notice that in the tests where a matter-field is present, we convert the conservative variables given in Eq.~\eqref{eqn.cons.var} back to the primitive variables $(\rho, v_{i}, p)$ following closely the procedure detailed in~\cite{DumbserZanottiGaburroPeshkov2023}.
For the test cases involving black hole spacetimes in horizon-penetrating coordinates, we employ the \textit{excision} technique to remove the singularity from the computational domain, adopting those excision box sizes that led to long-time stable simulations.
Regarding the boundary conditions at excision boundaries, we prescribe the exact solution in the interior of the excision box, so that we can apply the standard schemes in an unaltered fashion.
Finally, concerning the damping parameters that are inherent in the GH formulation, we have decided to choose test-case-dependent constants. In all cases, we have tried to use the least possible amount of damping, starting always with zero values as first choice.

\subsection{The robust stability test}
\label{sec:stability-test}
\begin{figure}[!htbp]
\begin{center}
    \includegraphics[width=0.9\textwidth, keepaspectratio]{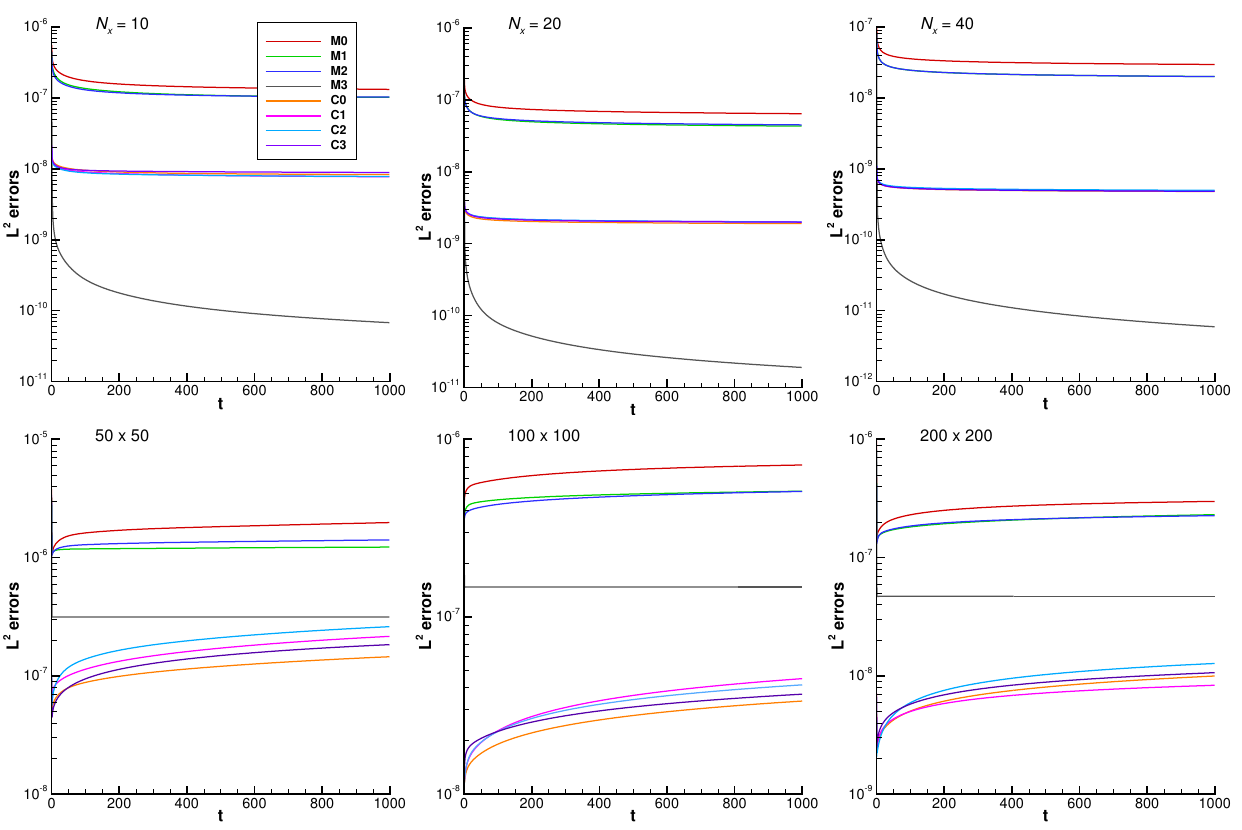}
	\caption{Time evolution of the constraints for the robust stability test case with a random initial perturbation of amplitude $\pm 10^{-7}/\varrho^2$ applied to all variables, performed on a sequence of successively refined meshes on the unit square in 2D. Top panel from left to right: fourth order ADER-DG scheme with $10$ ($\varrho = 1$), $20$ ($\varrho = 2$), and $40$ ($\varrho = 4$) elements in the $x$ direction. Bottom panel from left to right: fifth order CWENO-FD scheme on $50 \times 50$ ($\varrho = 1$), $100\times100$ ($\varrho = 2$), and $200\times200$ ($\varrho = 4$) grid points.}
	\label{fig.RobustStability}
\end{center}	
\end{figure}

A first validation procedure for numerical relativity codes is the so called {\emph{robust stability test}}, as described in~\cite{Alcubierre2004}. This test is designed to identify potential exponentially growing modes that could otherwise remain hidden for an extended evolution time, and therefore it is used to empirically verify the hyperbolicity of the governing system of equations. 
It is performed in a flat Minkowski spacetime without matter and it consists of applying small independent random perturbations to each evolved quantity of the PDE system. The amplitude of the perturbation is $\pm 10^{-7}/\varrho^2$, where the parameter $\varrho \in \left\{ 1, 2, 4 \right\}$ serves as mesh refinement factor. 
We underline that the chosen amplitude is three orders of magnitude larger than the one reported in~\cite{Alcubierre2004}, producing in this way an even more challenging configuration.
The computational domain is given by the unit square $\Omega = [-0.5, 0.5]\times[-0.5, 0.5]$, which is discretized setting $N_{x} = 10\varrho$ for the ADER-DG schemes, and by $50\varrho \times 50\varrho$  grid points for the CWENO-FD schemes. At the boundary, we impose flat Minkowski spacetime.
In Figure~\ref{fig.RobustStability} we show the evolution of the constraints for two sets of simulations, where we employ a fourth order ADER-DG and a fifth order CWENO-FD scheme. 
Here, no damping is used, i.e., we choose $\gamma_{0} = \gamma_{1} = \gamma_{2} = 0$. Moreover, we evolve the data in {\emph{harmonic gauge}}, setting $H_{a} = 0$ everywhere.
For a complete analysis, it is interesting to compare these results with those obtained using other first-order hyperbolic versions of well-known formulations, such as the Z4 system presented in~\cite{DumbserZanottiGaburroPeshkov2023}, and the BSSNOK system detailed in~\cite{DumbserZanottiPuppo2025}.

\subsection{The linearized gravitational wave}
\label{sec:linearized-wave}
The {\emph{linearized gravitational wave}} is essentially a one-dimensional test, consisting of a small gravitational wave perturbation of the flat Minkowski spacetime without matter, for which the metric results in
\begin{equation} 
	\label{eqn.lw.metric}
	d s^2 = - d t^2 + d x^2 + (1+b)\,d y^2 + (1-b)\,d z^2, \quad \textnormal{with} \quad b = \epsilon \sin \left( 2 \pi (x-t) \right)\,,
\end{equation} 
where $\epsilon = 10^{-8}$. This choice ensures that the model dynamics is linear and the terms depending on $\epsilon^{2}$ can be neglected.
As described in~\cite{Alcubierre2004}, this test is used to check the ability of a code to propagate a travelling gravitational wave, revealing possible sources of inaccuracy in the algorithms adopted. 

The components of the spatial metric can be readily identified from Eq.~\eqref{eqn.lw.metric} 
as $\gamma_{11}=1$, $\gamma_{22}=1+b$, and $\gamma_{33}=1{-}b$, as well as the ones of the shift vector, $\beta^{i} = 0$, and the value of the lapse function, $\alpha = 1$. Other relevant non-zero quantities 
worth mentioning are $\Pi_{22}$ and $\Pi_{33}$, due to their direct relation with the extrinsic curvature $K_{ij}$, which in the no-shift case reduces to $\Pi_{ij} = 2 K_{ij}$, see Eq.~\eqref{eq:Pi-K}.
We consider a rectangular domain $\Omega = [-0.5,0.5]\times[-0.05,0.05]$ with periodic boundary conditions in both directions, and we use the {\emph{harmonic gauge condition}} $H_{a} = 0$, while all damping parameters are set to zero.

\begin{figure}[!htbp]
\begin{center}
    \includegraphics[width=0.8\textwidth, keepaspectratio]{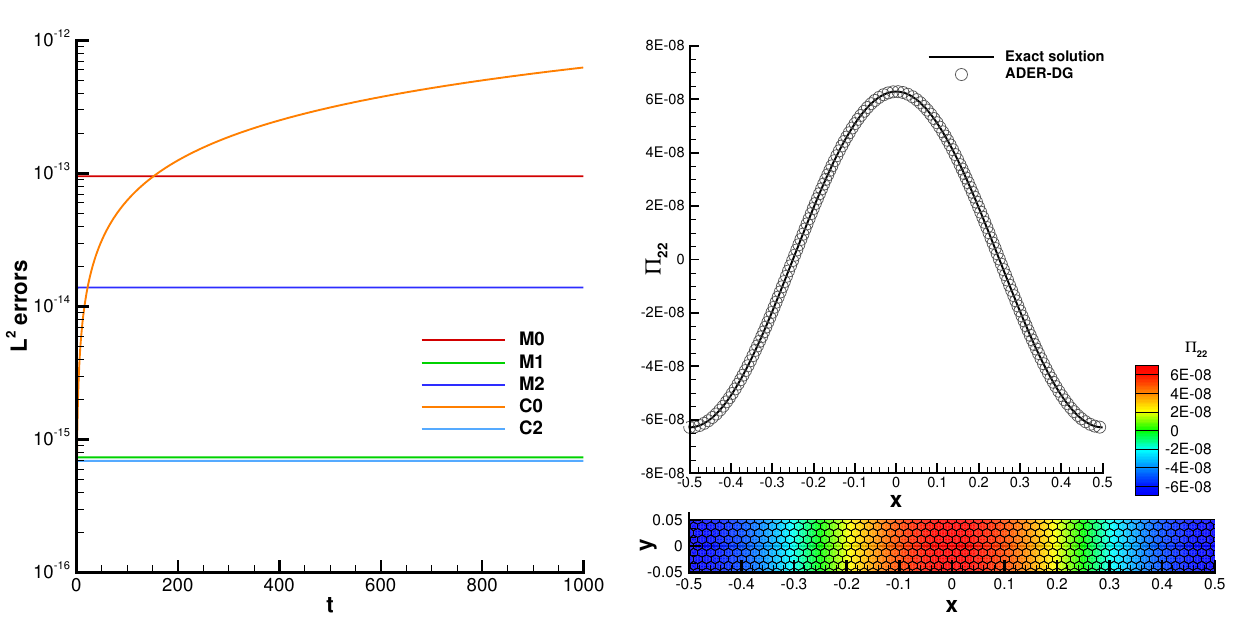}
	\caption{Solution of the linearized wave test using the fourth order ADER-DG scheme on an unstructured polygonal mesh. Left panel: time evolution of the constraints. Right panel: below, 2D color contour of $\Pi_{22}$; above, comparison of the exact and numerical solution at final time $t = 1000$ on a one dimensional cut at $y = 0$.}
	\label{fig.LinearizedWave.AleVoronoi.P3P3}
\end{center}	
\end{figure}

In Figure~\ref{fig.LinearizedWave.AleVoronoi.P3P3} we present the results obtained with the ADER-DG scheme of order four on a mesh where $N_{x} = 60$. In particular, we monitor the time evolution of the constraints, and we report the profile of the variable $\Pi_{22}$ over the entire 2D domain $\Omega$, together with a one-dimensional cut that is compared against the exact solution.

Since the CWENO-FD schemes provide results similar to those of Figure~\ref{fig.LinearizedWave.AleVoronoi.P3P3}, we omit the corresponding plots. Instead,  in Table~\ref{tab.LinearizedWave.FDToy.convergence}, we perform a numerical convergence analysis at short time $t = 1$ to show that they achieve the expected order of accuracy.  

In all the following convergence tables, we will use the notation $\mathbb{P}_{N}$ to denote, in the CWENO-FD case, the schemes obtained by reconstructing an optimal polynomial of degree $N$, see Eq.~\eqref{eq.Popt}, while in the ADER-DG case, the schemes obtained by representing the discrete solution by piecewise polynomials of degree $N$, see Eq.~\eqref{eqn.uh}.

\begin{table}[!htbp]
\centering
\begin{tabular}{lcccc}
\toprule
\textbf{Scheme} & $\bm{N_{x}\times N_{y}}$ & $\bm{L^2}$ \textbf{Error} $\bm{\Pi_{22}}$ & $\bm{\mathcal{O}(\Pi_{22})}$ & \textbf{Theoretical}\\
\midrule
CWENO-FD $\mathbb{P}_{4}$ & $8\times8$ & 3.9347E-10 & - \\
 & $16\times16$ & 1.3403E-11 & 4.88 & 5.00\\
 & $24\times24$ & 1.7898E-12  & 4.97 \\
\midrule
CWENO-FD $\mathbb{P}_{6}$ & $8\times8$ & 5.0194E-11 & -  \\
 & $16\times16$ & 4.3791E-13 & 6.84 & 7.00 \\
 & $24\times24$ & 2.6151E-14 &  6.95 \\
\midrule
CWENO-FD $\mathbb{P}_{8}$ & $10\times10$ & 9.4470E-13 & - \\
 & $16\times16$ & 1.4826E-14 & 8.84 & 9.00\\
 & $24\times24$ & 3.9615E-16 & 8.93 \\
\bottomrule
\end{tabular}
\caption{Numerical convergence results for the linearized wave test at $t=1$ for the variable $\Pi_{22}$ using the CWENO-FD schemes.
\label{tab.LinearizedWave.FDToy.convergence}
}
\end{table}

\subsection{The gauge wave}
\label{sec:gauge-wave}
Another significant test taken from~\cite{Alcubierre2004} is the so called {\emph{gauge wave test}}, which involves the evolution of a propagating wave on a flat Minkowski spacetime with a time-dependent spatial metric. A common choice, which we adopt here, is to consider a sinusoidal profile of amplitude $A$ propagating along the $x$-axis, leading to
\begin{equation}
	\label{eqn.gw.metric}
	d s^2 = - H(x,t) \, d t^2 + H(x,t) \, d x^2 + d y^2 + d z^2, \quad \text{where} \quad H(x,t) = 
	1-A\,\sin \left( 2\pi(x-t) \right)\,.
\end{equation}
Consequently, we have $\gamma_{11}=H$ and $\gamma_{22}= \gamma_{33} = 1$, with no shift, i.e. $\beta^{i} = 0$, and a time-dependent lapse function $\alpha = \sqrt{H}$. The only non-zero components of the tensor $\Pi_{ab}$ are $\Pi_{00}$ and $\Pi_{11}$, where $\Pi_{00} = - \Pi_{11}$ with 
\begin{equation}
	\Pi_{11} = - 2\pi A\frac{\cos\left(2\pi(x-t)\right)}{\sqrt{1 - A\sin\left(2\pi(x-t)\right)}}\,.
\end{equation}
We run our tests over a rectangular domain $\Omega = [-0.5,0.5]\times[-0.05,0.05]$ with periodic boundary conditions in both directions. Moreover, we again use the {\emph{harmonic gauge condition}} $H_{a} = 0$, and no damping whatsoever, that is, $\gamma_{0} = \gamma_{1} = \gamma_{2} = 0$. 

At first glance, there seems to be no reason to consider this test particularly challenging. However, it is well-documented that its successful resolution is, in fact, quite involved. Although it depends on the specific formulation under analysis, difficulties typically become evident at high amplitude, i.e. $A = 0.1$.
This is the case for the BSSNOK system, which fails after a rather short time in both its first and second order classical formulations, as described in~\cite{Alic:2011a,Brown2012}. The original version of the CCZ4 system is stable only in its damped formulation~\cite{Alic:2011a}, whereas the first-order formulation FO-CCZ4 presented in~\cite{Dumbser2017strongly} provided the first stable undamped simulation with CCZ4. Conversely, this test can be effectively handled within the Z4 formulation~\cite{DumbserZanottiGaburroPeshkov2023} as well as the GH formulation~\cite{Babiuc2006, Boyle2007, lovelace2025simulating}, which is the object of study in this paper.

\begin{figure}[!htbp]
\begin{center}
    \includegraphics[width=0.8\textwidth, keepaspectratio]{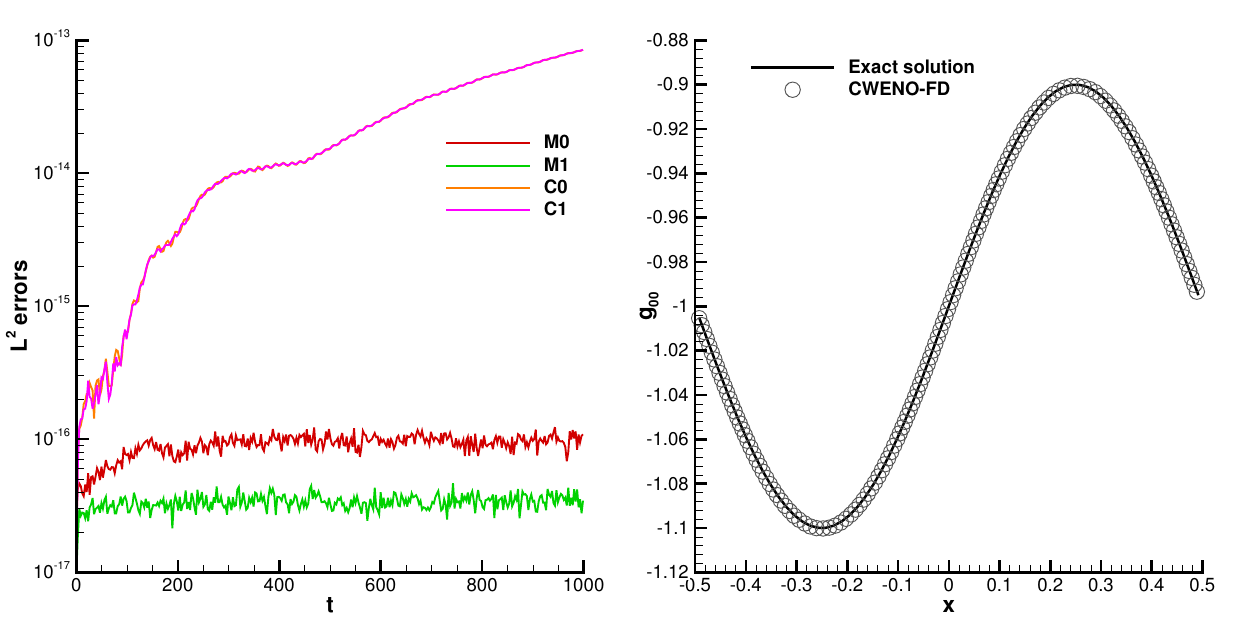}
	\caption{Solution of the gauge wave test with amplitude $A = 0.1$ using the fifth order CWENO-FD scheme on a Cartesian mesh. Left panel: time evolution of the constraints. Right panel: comparison of the exact and numerical solution for the variable $g_{00}$ at final time $t = 1000$ on a one dimensional cut at $y = 0$.}
	\label{fig.GaugeWave.CWENO5}
\end{center}	
\end{figure}

\begin{figure}[!htbp]
\begin{center}
    \includegraphics[width=0.8\textwidth, keepaspectratio]{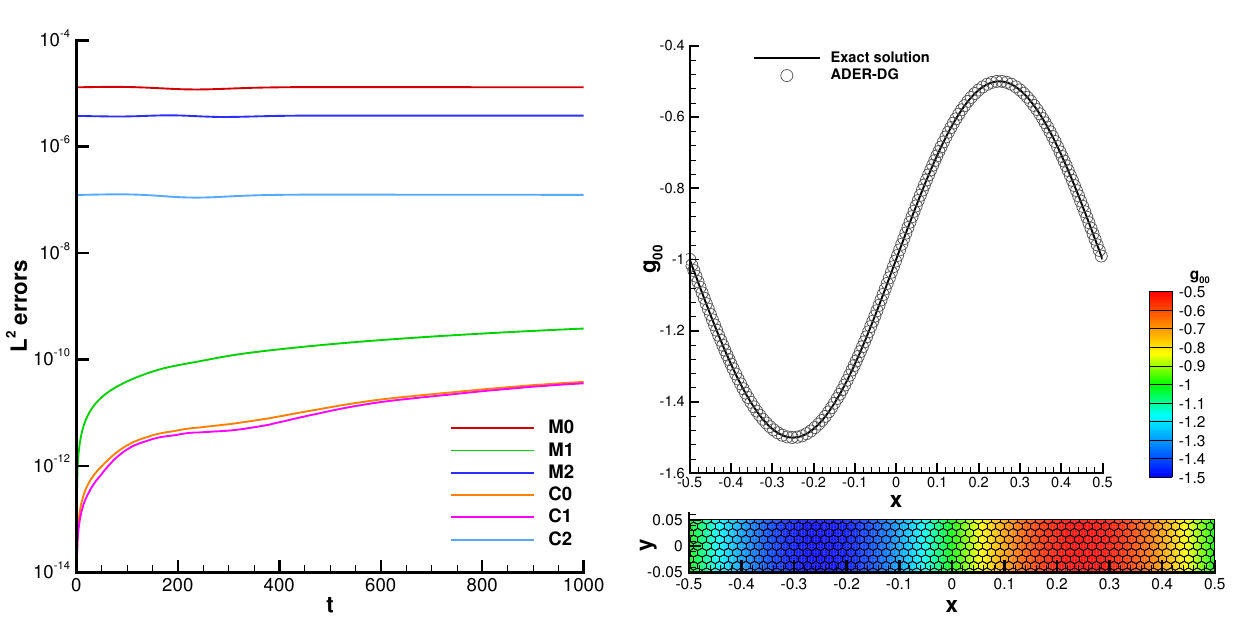}
	\caption{Solution of the gauge wave test with high amplitude $A = 0.5$ using the fourth order ADER-DG scheme on an unstructured polygonal mesh. Left panel: time evolution of the constraints. Right panel: below, 2D color contour of $g_{00}$; above, comparison of the exact and numerical solution at final time $t = 1000$ on a one dimensional cut at $y = 0$.}
	\label{fig.GaugeWave.AleVoronoi.P3P3}
\end{center}	
\end{figure}

We first perform a simulation with large amplitude $A=0.1$ using the fifth order CWENO-FD scheme on a $64 \times 16$ grid, with corresponding results  presented in Figure~\ref{fig.GaugeWave.CWENO5}. We then run a second simulation with an even higher amplitude $A=0.5$ employing the fourth order ADER-DG method on a mesh generated by $N_{x} = 60$ elements in the $x$ direction, as shown in Figure~\ref{fig.GaugeWave.AleVoronoi.P3P3}. In both cases, 
we notice that the constraints remain controlled throughout the time evolution, and the exact profile of the variable $g_{00}$ is perfectly recovered at the final time $t = 1000$, suggesting that both tests are well resolved.

To conclude our study, in Tables~\ref{tab:GaugeWave.FDToy.convergence} and \ref{tab:GaugeWave.AleVoronoi.convergence} we present a numerical convergence analysis at time $t = 1$ and high amplitude $A = 0.5$ for the CWENO-FD and ADER-DG schemes, respectively. In the latter case, we also provide the maximum diameter $h$ of the circumcircles over all elements as an additional characteristic measure of the mesh size.

\begin{table}[!htbp]
\centering
\begin{tabular}{lccccc}
\toprule
\textbf{Scheme} & $\bm{N_{r}\times N_{\theta}}$ & $\bm{L^2}$ \textbf{Error} $\bm{g_{00}}$ & $\bm{\mathcal{O}(g_{00})}$ & \textbf{Theoretical}\\
\midrule
CWENO-FD $\mathbb{P}_{4}$  & $32\times32$ &  1.3335E-05 & - & \\
 & $64\times64$ & 4.5713E-07 & 4.87 & 5.00\\
 & $96\times96$ & 6.1146E-08 & 4.96 & \\
\midrule
CWENO-FD $\mathbb{P}_{6}$ & $32\times32$ & 1.6025E-06 & -  \\
 & $64\times64$ & 1.5767E-08 & 6.67 & 7.00 \\
 & $96\times96$ & 9.6039E-10 & 6.90 & \\
 \midrule
CWENO-FD $\mathbb{P}_{8}$ & $32\times32$ & 3.4042E-07 & -\\
 & $64\times64$ & 1.0864E-09 & 8.29 & 9.00\\
 & $96\times96$ & 3.8951E-11 & 8.21 &  \\
\bottomrule
\end{tabular}
\caption{Numerical convergence results for the gauge wave test with high amplitude $A = 0.5$ at $t=1$ for the variable $g_{00}$ using the family of CWENO-FD schemes on Cartesian grids.
}
\label{tab:GaugeWave.FDToy.convergence}
\end{table}

\begin{table}[!htbp]
\centering
\begin{tabular}{lccccc}
\toprule
\textbf{Scheme} & $\bm{N_{x}}$ & $\bm{h}$ & $\bm{L^2}$ \textbf{Error} $\bm{g_{00}}$ & $\bm{\mathcal{O}(g_{00})}$ & \textbf{Theoretical}\\
\midrule
ADER-DG $\mathbb{P}_{1}$ & $40$ & 8.3140E-03 & 1.7907E-04 & - \\
 & $80$ & 4.4844E-03 & 4.4628E-05 &  2.25 & 2.00\\
 & $120$ & 3.0898E-03 & 2.0248E-05 & 2.12\\
\midrule
ADER-DG $\mathbb{P}_{2}$ & $40$ & 8.3140E-03  & 2.0386E-06 & -  \\
 & $80$ & 4.4844E-03  & 2.6027E-07 & 3.33 & 3.00\\
 & $120$ & 3.0898E-03 & 7.7873E-08 & 3.24 \\
\midrule
ADER-DG $\mathbb{P}_{3}$ & $40$ & 8.3140E-03 & 1.5630E-08 & - \\
 & $80$ & 4.4844E-03 & 1.0416E-09 & 4.39 & 4.00\\
 & $120$ & 3.0898E-03 & 2.0965E-10 & 4.30 \\
 \midrule
ADER-DG $\mathbb{P}_{4}$ & $20$ & 2.0127E-02 & 5.5749E-09 & - \\
 & $40$ & 8.3140E-03  & 1.7908E-10 & 3.89 & 5.00\\
 & $80$ & 4.4844E-03 & 9.9518E-12 & 4.68 \\
\bottomrule
\end{tabular}
\caption{Numerical convergence results for the gauge wave test with high amplitude $A = 0.5$ at $t=1$ for the variable $g_{00}$ using the family of ADER-DG schemes on unstructured polygonal meshes.
}
\label{tab:GaugeWave.AleVoronoi.convergence}
\end{table}

\subsection{Single stationary black holes in two and three space dimensions}

The first exact nontrivial solution of the Einstein equations is the {\emph{Schwarzschild metric}}, derived by 
Karl Schwarzschild~\cite{1916SPAW.......189S} shortly after Einstein published his foundational paper on the theory of general relativity.
It describes the spacetime surrounding a static, spherically symmetric object, such as a non-rotating black hole.
Then, almost fifty years had to pass before a generalization to rotating, stationary, axially symmetric black holes was proposed by Roy Kerr~\cite{Kerr63}, reflecting the difficulties of its derivation. This solution is known as the {\emph{Kerr metric}}.

Throughout the simulations reported here, the black hole mass is chosen as $M=1\,M_{\odot}$, and $T_{ab}$ vanishes since we are indeed dealing with vacuum solutions. Moreover, we remark that the 2D tests are performed using the ADER-DG scheme, while for the 3D ones the CWENO-FD method is employed instead.

\begin{figure}[!htbp]
\begin{center}
    \includegraphics[width=0.8\textwidth, keepaspectratio]{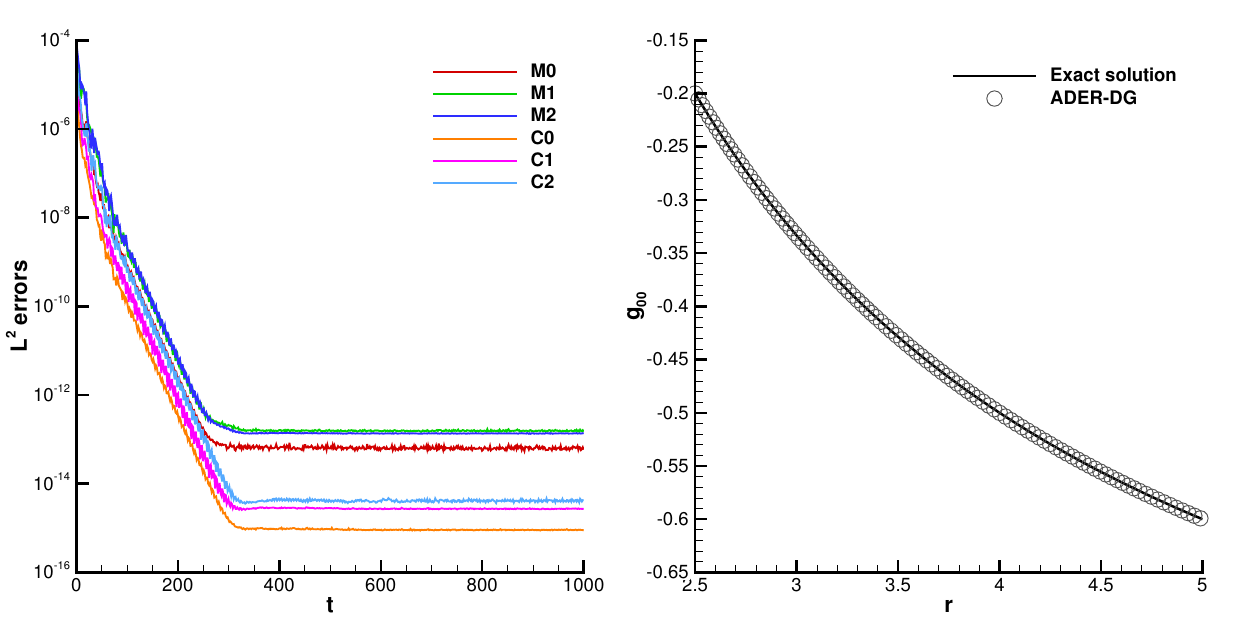}
	\caption{2D Schwarzschild black hole subject to an initial Gaussian perturbation in the variable $g_{00}$ using the well-balanced fourth order ADER-DG scheme. Left panel: time evolution of the constraints, which decay after a relatively short period of time. Right panel: comparison of the exact and numerical solution for the variable $g_{00}$ at final time $t = 1000$ on a one dimensional cut at $\theta = \pi/2$.}
	\label{fig.SchwarzschildBL.WB.AleVoronoi.P3P3}
\end{center}	
\end{figure}

\paragraph{Non-rotating black hole in 2D}
We begin by evolving a Schwarzschild black hole ($a=0$) in Boyer-Lindquist coordinates (see, for instance, Appendix B of~\cite{Komissarov04b}, or the more extended discussion on the different coordinate systems adapted to the Kerr spacetime in~\cite{Visser2007}).
The two-dimensional computational domain $(r, \theta) \in [2.5, 5] \times [\delta,\pi-\delta]$, with $\delta = 0.520796327$, is discretized by $50$ elements in the $r$ direction, namely setting $N_{x} = 50$. At both boundaries, we impose Dirichlet boundary condition for all variables, using the initial configuration. The gauge source function $H_{a}$ is computed directly from the analytical initial data as $H_{a} = -\Gamma_{a}$, and we set the damping parameters to $\gamma_{0} = \gamma_{2} = 1, \gamma_{1} = -1$.
We introduce a {\emph{small Gaussian perturbation}} on the $g_{00}$ component of the initial spacetime metric with the aim of studying the well-balancing property of the scheme. More precisely, at time $t = 0$, we set
\begin{equation}\label{eq.gaussian.pert}
	\tilde{g}_{00}(r,\theta) = g_{00}(r,\theta) + A \exp \left( -\frac{1}{2} \frac{(r\sin\theta-4)^2 + (r\cos\theta-0)^2}{\sigma^2} \right)\,,
\end{equation}   
with $A = 10^{-4}$ and $\sigma=0.2$.
We present the results obtained with the well-balanced fourth order version of the ADER-DG scheme in Figure~\ref{fig.SchwarzschildBL.WB.AleVoronoi.P3P3}. In line with our expectations, the perturbation gradually propagates out of the domain, so that at sufficiently large time the solution returns to the exact stationary equilibrium state.

We then provide a further setup, in which we fully exploit the ability of our ADER-DG code to operate on general unstructured polygonal meshes. To this end, we evolve a Schwarzschild black hole using the 3D Cartesian horizon-penetrating harmonic coordinates $(t,x,y,z)$ described in SpECTRE documentation \cite{deppe_2025_16906840}, which are in turn based on \cite{cook1997well}. To perform a 2D simulation, however, we restrict our attention to the evolution of the solution over the equatorial plane $z = 0$. We emphasize that, in this configuration, we need to take special care of the $z$ direction by introducing an artificial source term in Eq.~\eqref{eqn.pde.mat.preview} to account for possibly non-zero terms arising from derivatives of the state variables with respect to the $z$ coordinate. 
More precisely, we \textit{modify} the system~\eqref{eqn.pde.mat.preview} as follows
\begin{equation}
	\label{eqn.pde.mod}
	\frac{\partial \u }{\partial t} + \nabla \cdot \F(\u)
	+ \B(\u) \cdot \nabla \u  
	= \S(\u) + \S^a(\x),
\end{equation}
where the artificial source term $\S^a(\x)$ is based just on the knowledge of the exact equilibrium solution $\u_e(\x)$ and is given by 
\begin{equation}
	\label{eqn.artificial.source}
	\S^a(\x) = - \partial_z \mathbf{f}^z(\u_e(\x)) - \mathbf{B}^z(\u_e(\x)) \partial_z \u_e(\x). 
\end{equation}
Such term is evaluated at $z = 0$ as we focus on the equatorial plane, and it accounts for all spatial derivatives in the $z$ direction that are not computed in a 2D code. 
This passage is in fact crucial, otherwise the intrinsically 3D nature of the coordinate system would be inconsistent with the 2D nature of our code, and thus the Einstein equations would not be completely satisfied due to the omission of these contributions; see also \cite{Gaburro2021WBGR1D} for a similar approach in 1D. Notice that adding this artificial source term does not alter the structure of the principal part of the PDE, i.e. hyperbolicity remains guaranteed.

\begin{figure}[!htbp]
\begin{center}
    \includegraphics[width=0.85\textwidth, keepaspectratio]{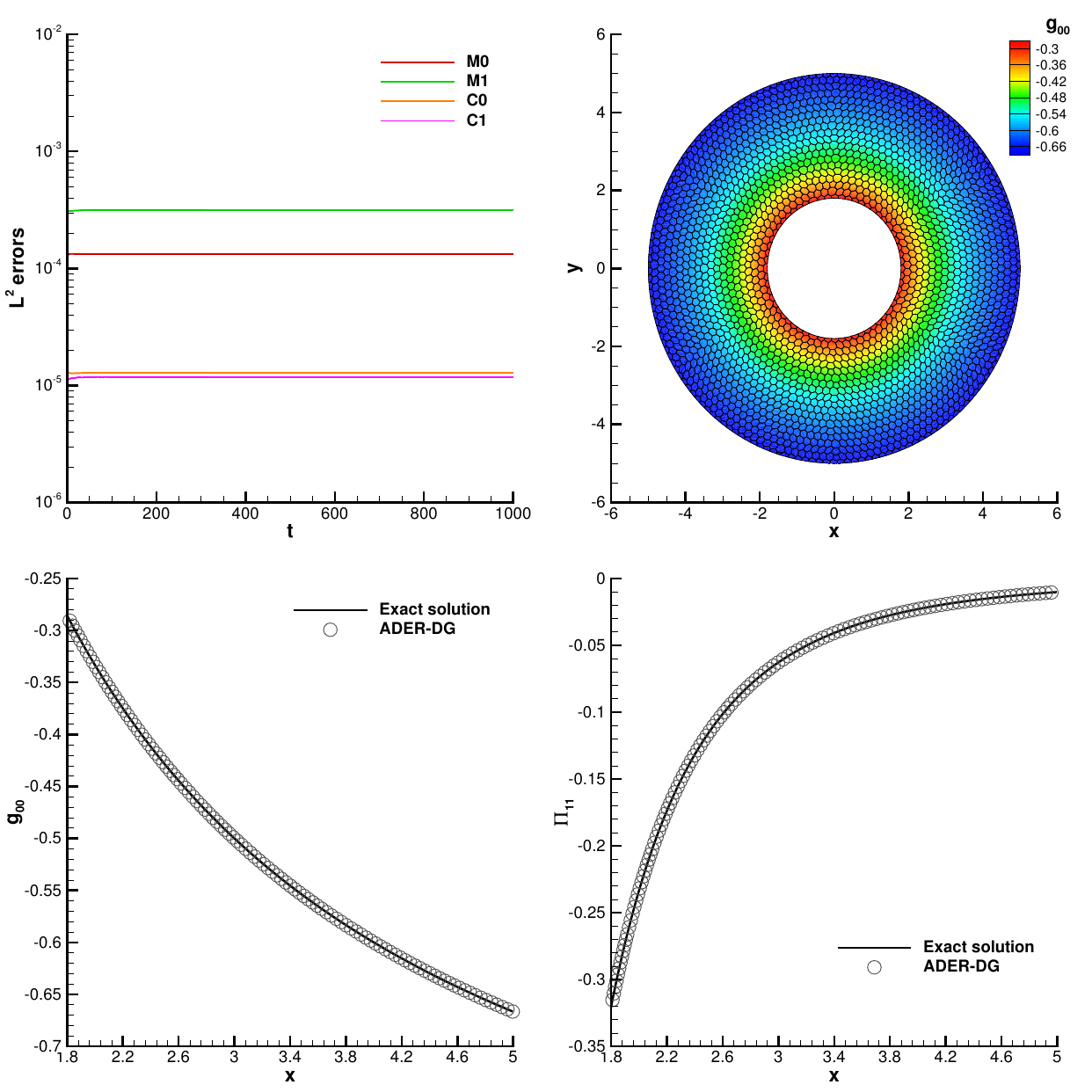}
	\caption{2D Schwarzschild black hole in 3D Cartesian harmonic coordinates evolved over the equatorial plane $z = 0$ using the fourth order ADER-DG scheme. Top-left panel: time evolution of the constraints. Top-right panel: 2D color contour of $g_{00}$  at final time $t = 1000$, including a detail of the polygonal mesh elements. Bottom panels: comparison of the exact and numerical solution for the variable $g_{00}$ and $\Pi_{11}$ at final time $t = 1000$ on a one dimensional cut at $y = 0$.}
	\label{fig.Schwarzschild.disk.AleVoronoi.P3P3}
\end{center}	
\end{figure}

As computational domain, we consider an annulus $\Omega$ with inner radius $r_{\min} = 1.8$ and external radius $r_{\max} = 5.0$, discretized by a total number of elements equal to $N_{P} = 2257$, with characteristic mesh size of $h = \frac{|\Omega|}{N_{P}} = \frac{\pi (r_{\max}^{2} -\, r_{\min}^{2})}{N_{P}} = 3.0288 \times10^{-2}$. At both boundaries, we impose the analytical initial data. Notice also how it is the mesh geometry itself that incorporates the excision zone, which is simply determined by the choice of $r_{\min}$. Since the coordinates are harmonic in time and space, we use the \emph{harmonic gauge condition} for $H_{a}$. Moreover, we set the damping parameters to $\gamma_0 = \gamma_2 = 1, \gamma_1 = -1$. Figure \ref{fig.Schwarzschild.disk.AleVoronoi.P3P3} shows the results of this test case obtained using the fourth order version of the ADER-DG schemes, together with a visual representation of the computational domain $\Omega$ and the associated unstructured polygonal mesh.

\begin{figure}[!htbp]
\begin{center}
    \includegraphics[width=0.8\textwidth, keepaspectratio]{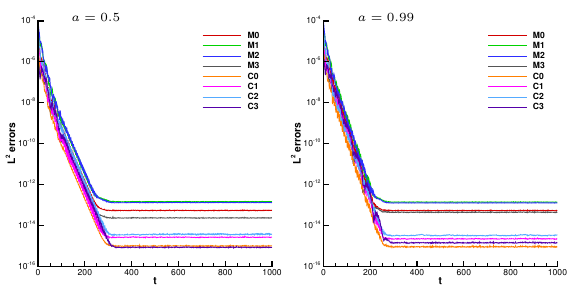}
	\caption{2D Kerr black holes subject to an initial Gaussian perturbation in the variable $g_{00}$ using the well-balanced fourth order ADER-DG scheme. Time evolution of the constraints for spin $a = 0.5$ on the left, and for extreme spin $a = 0.99$ on the right. In both cases the perturbation decays so that the exact equilibrium state is recovered.}
	\label{fig.KerrBL.05.099.WB.AleVoronoi.P3P3}
\end{center}	
\end{figure}

\paragraph{Rotating black hole in 2D}
Next, we evolve two Kerr black holes, one with spin $a = 0.5$, and the other one with extreme spin $a = 0.99$, using the same coordinates, setting, and resolution as in the first 2D non-rotating test case described above. We also retain the same structure for the numerical test, that is, we apply a Gaussian perturbation to the variable $g_{00}$ as in Eq.~\eqref{eq.gaussian.pert}. 
Figure~\ref{fig.KerrBL.05.099.WB.AleVoronoi.P3P3} shows the time evolution of the constraints for the two cases, confirming again that the well-balancing feature allows the solutions to come back again to the correct equilibrium states up to machine precision.

\begin{figure}[!htbp]
\begin{center}
    \includegraphics[width=\textwidth, keepaspectratio]{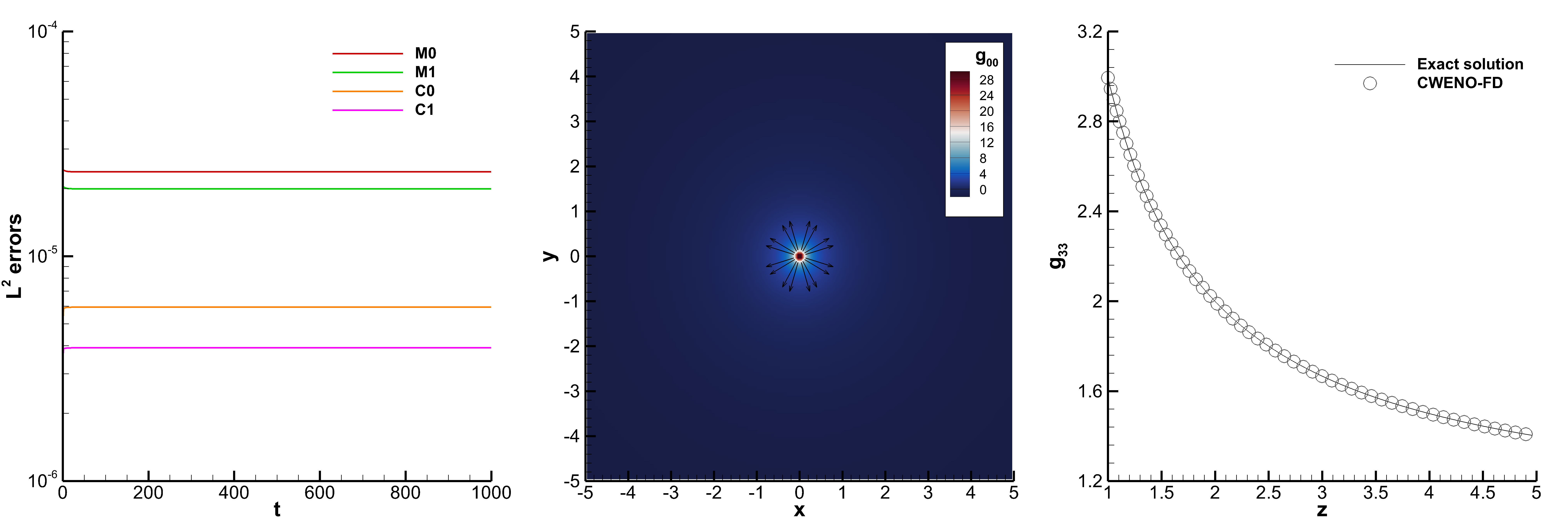}
	\caption{3D Schwarzschild black hole in Cartesian Kerr-Schild coordinates using the seventh order CWENO-FD scheme. Left panel: time evolution of the constraints. Central panel: 2D cut on the equatorial plane $z = 0$ at final time $t = 1000$ with detail of the shift vector field near the singularity. Right panel: comparison of the exact and numerical solution for the variable $g_{33}$ at final time $t = 1000$ on a one dimensional cut along the $z$-axis.}
	\label{fig.Schwarzschild3D.noWB.FDToy.CWENO7}
\end{center}	
\end{figure}

\paragraph{Non-rotating black hole in 3D}
We now turn to tests in three space dimensions using the 3D Cartesian Kerr-Schild coordinates~\cite{Visser2007}, first by evolving a stationary Schwarzschild black hole ($a = 0$).
Here the computational domain is defined by $\Omega = [-5, 5]^{3}$, from which we have excised a cubic volume with a half-edge length of unity, centered on the physical singularity at $r=0$, i.e. $\Omega_e = [-1, 1]^{3}$.
Figure~\ref{fig.Schwarzschild3D.noWB.FDToy.CWENO7} shows the results obtained using the seventh order version of our CWENO-FD schemes with a resolution of $132^3$ grid points, and damping parameters set to $\gamma_{0} = \gamma_{2} = 1, \gamma_{1} = -1$. Notice that the central panel, which represents a 2D cut on the equatorial plane $z = 0$ with detail of the shift vector field near the singularity, is mainly used as a visual representation, as the cube $\Omega_e$ is not evolved, corresponding indeed to the excision region.

\begin{figure}[!htbp]
\begin{center}
    \includegraphics[width=\textwidth, keepaspectratio]{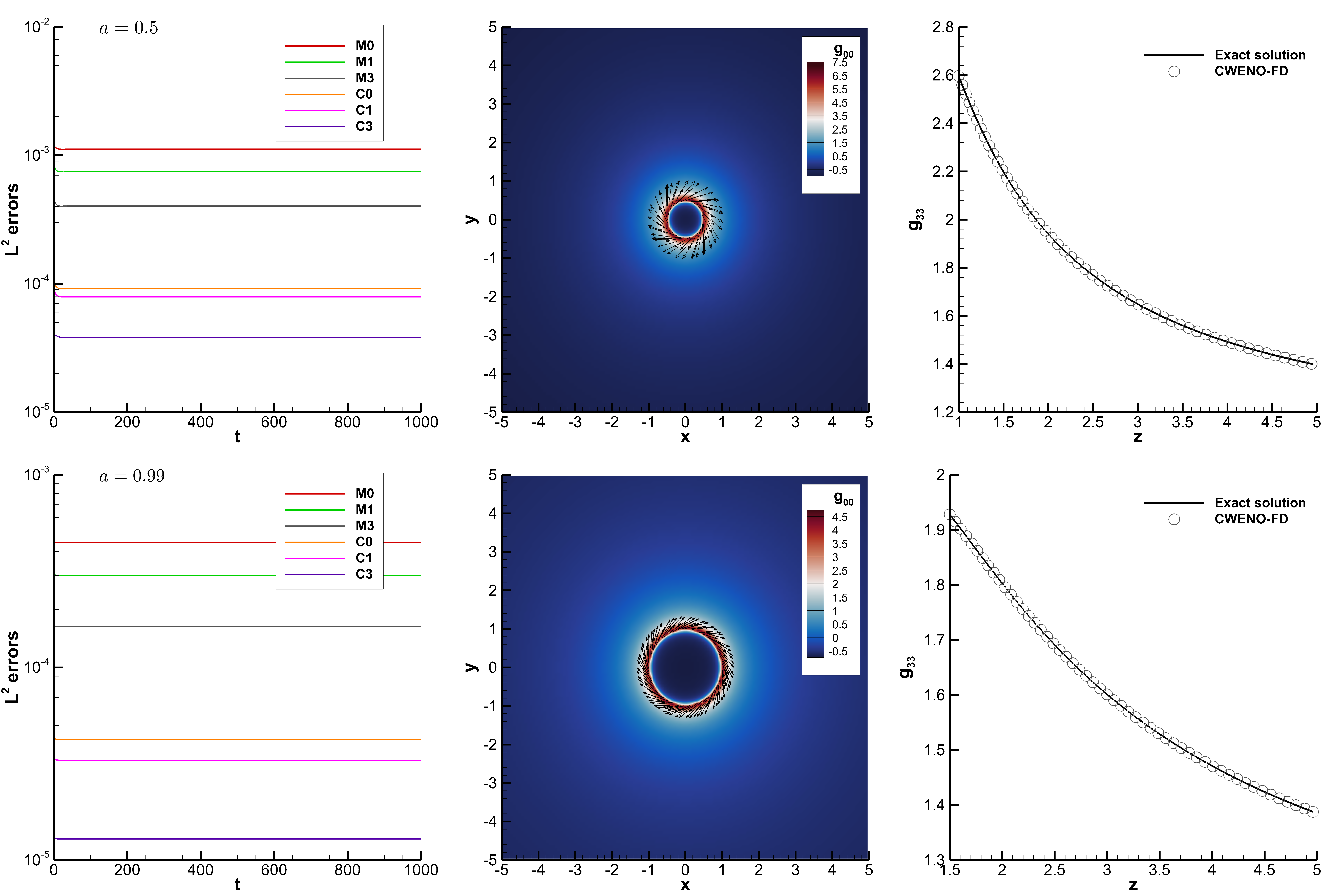}
	\caption{3D Kerr black hole in Cartesian Kerr-Schild coordinates using the seventh order CWENO-FD scheme. Top panels: spin $a = 0.5$. Bottom panels: spin $a = 0.99$. Left panels: time evolutions of the constraints. Central panels: 2D cuts on the equatorial plane $z = 0$ at final time $t = 1000$ with detail of the shift vector fields near the ring singularities. Right panels: comparisons of the exact and numerical solutions for the variable $g_{33}$ at final time $t = 1000$ on one dimensional cuts along the $z$-axis.
}
\label{fig.Kerr3D.05.099.noWB.FDToy.CWENO7.combined}
\end{center}	
\end{figure}

\paragraph{Rotating black hole in 3D}
To conclude, we consider two Kerr black holes in three space dimensions, with spins $a = 0.5$ and $a = 0.99$. The domain, resolution, damping parameters, and scheme adopted are the same as in the 3D Schwarzschild case. Notice that, unlike the point singularity in the Schwarzschild metric, the Kerr singularity forms a ring with radius $a$ lying in the $z = 0$ plane, the so called {\emph{ring singularity}}, which we must enclose within the excision region. To this end, we excise cubic volumes with half-edge lengths of $1$ and $1.5$ for $a = 0.5$ and $a = 0.99$, respectively, both centered at the origin.
In Figure~\ref{fig.Kerr3D.05.099.noWB.FDToy.CWENO7.combined} we provide the results associated to the two spins discussed.

\subsection{Special relativistic Riemann problems in Cowling approximation}
\label{sec:Riemann-problems}
To evaluate the robustness of our numerical schemes in handling the matter evolution, we consider a series of Riemann problems within the Cowling approximation~\cite{Cowling41}, which consists on neglecting the evolution of spacetime, here represented by a flat Minkowski metric in Cartesian coordinates. This approach provides a controlled framework to isolate the performance of the relativistic hydrodynamic solvers, particularly in resolving strong shocks and discontinuities, without the additional complexity of spacetime dynamics.
We therefore address three Riemann problems, previously studied for instance by~\cite{DumbserZanottiPeshkov2024,DumbserZanottiPuppo2025}, with the following initial conditions and resulting wave structures:
\begin{itemize}
\item $(\rho, \, v_{x}, \, p)_L = (1,\, -0.6, \,10)$ and $(\rho, \, v_{x}, \, p)_R = (10, \, 0.5, \,20)$ with adiabatic index $\gamma = 5/3$, as considered by~\cite{Mignone2005}. We refer to this as the \textit{2R problem}, since it produces two rarefaction waves;

\item $(\rho, \, v_{x}, \, p)_L = (10^{-3}, \, 0, \, 1)$ and $(\rho, \, v_{x}, \, p)_R = (10^{-3}, \, 0, \, 10^{-5})$ with adiabatic index $\gamma = 5/3$, as considered by~\cite{Radice2012a}. We refer to this as the  \textit{RS problem}, since it produces a rarefaction wave and a shock wave;

\item $(\rho, \, v_{x}, \, p)_L = (1, \,0.9, \,1)$ and $(\rho, \, v_{x}, \, p)_R = (1, \,0, \,10)$ with adiabatic index $\gamma = 4/3$, as considered by~\cite{Zhang2006}. We refer to this as the  \textit{2S problem}, since it produces two shock waves.
\end{itemize}

\begin{figure}[!htbp]
	\begin{center}
		\includegraphics[width=\textwidth, keepaspectratio]{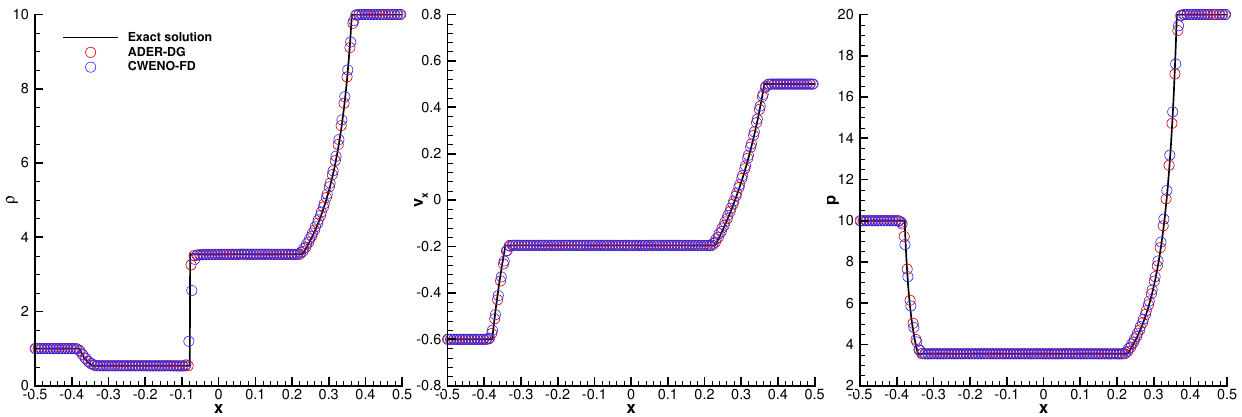}
		\caption{
		Solution of the \textit{2R problem} at final time $t = 0.4$ using the third order ADER-DG and CWENO-FD schemes. Left, central, and right panels: comparison of the exact and numerical solutions for the primitive variables $\rho, v_{x}$, and $p$, respectively, along a one dimensional cut at $y = 0$.
		}%
		\label{fig.RHD.Riemann1.2R}
	\end{center}	
\end{figure}

\begin{figure}[!htbp]
	\begin{center}
		\includegraphics[width=\textwidth, keepaspectratio]{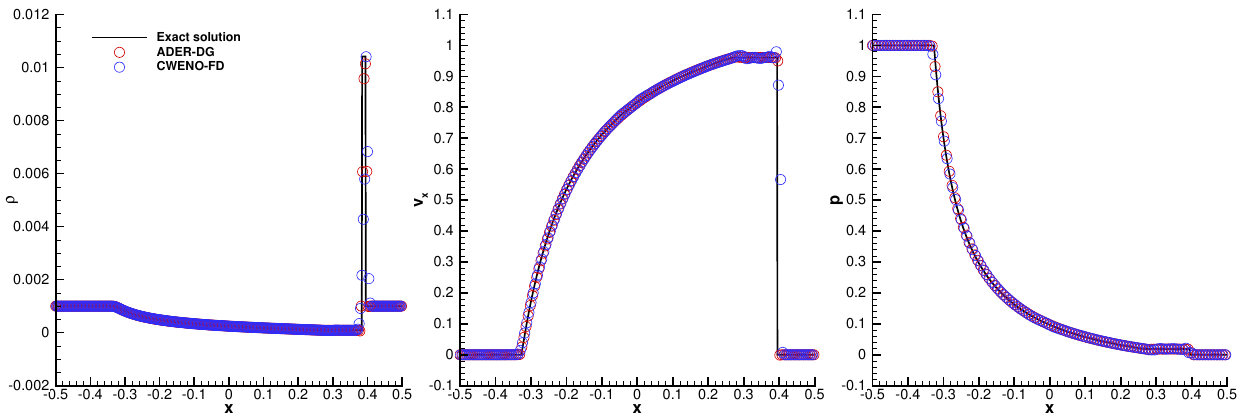}
		\caption{
		Solution of the \textit{RS problem} at final time $t = 0.4$ using the third order ADER-DG and CWENO-FD schemes. Left, central, and right panels: comparison of the exact and numerical solutions for the primitive variables $\rho, v_{x}$, and $p$, respectively, along a one dimensional cut at $y = 0$.
		}%
		\label{fig.RHD.Riemann2.RS}
	\end{center}	
\end{figure}

\begin{figure}[!htbp]
	\begin{center}
		\includegraphics[width=\textwidth, keepaspectratio]{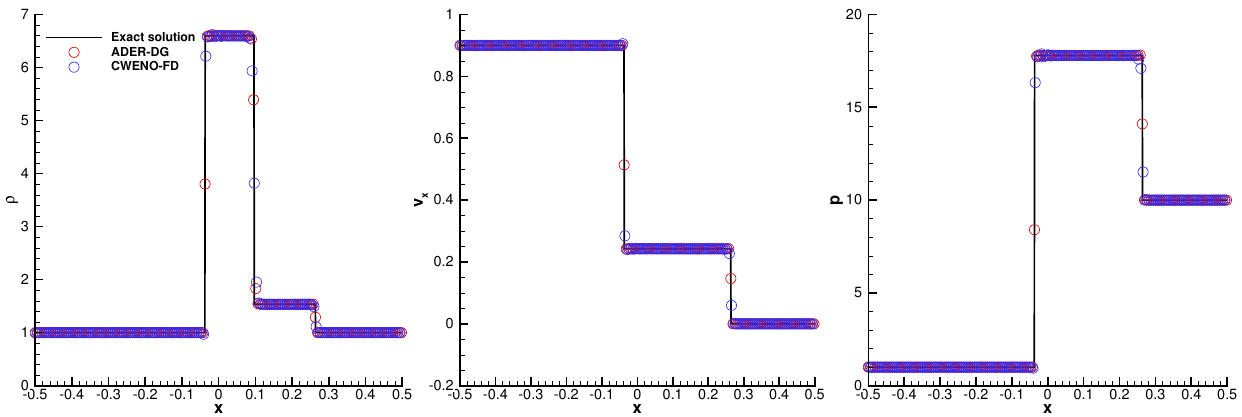}
		\caption{
		Solution of the \textit{2S problem} at final time $t = 0.4$ using the third order ADER-DG and CWENO-FD schemes. Left, central, and right panels: comparison of the exact and numerical solutions for the primitive variables $\rho, v_{x}$, and $p$, respectively, along a one dimensional cut at $y = 0$.
		}%
		\label{fig.RHD.Riemann3.2S}
	\end{center}	
\end{figure}

\begin{figure}[!htbp]
	\begin{center}
		\includegraphics[width=\textwidth, keepaspectratio]{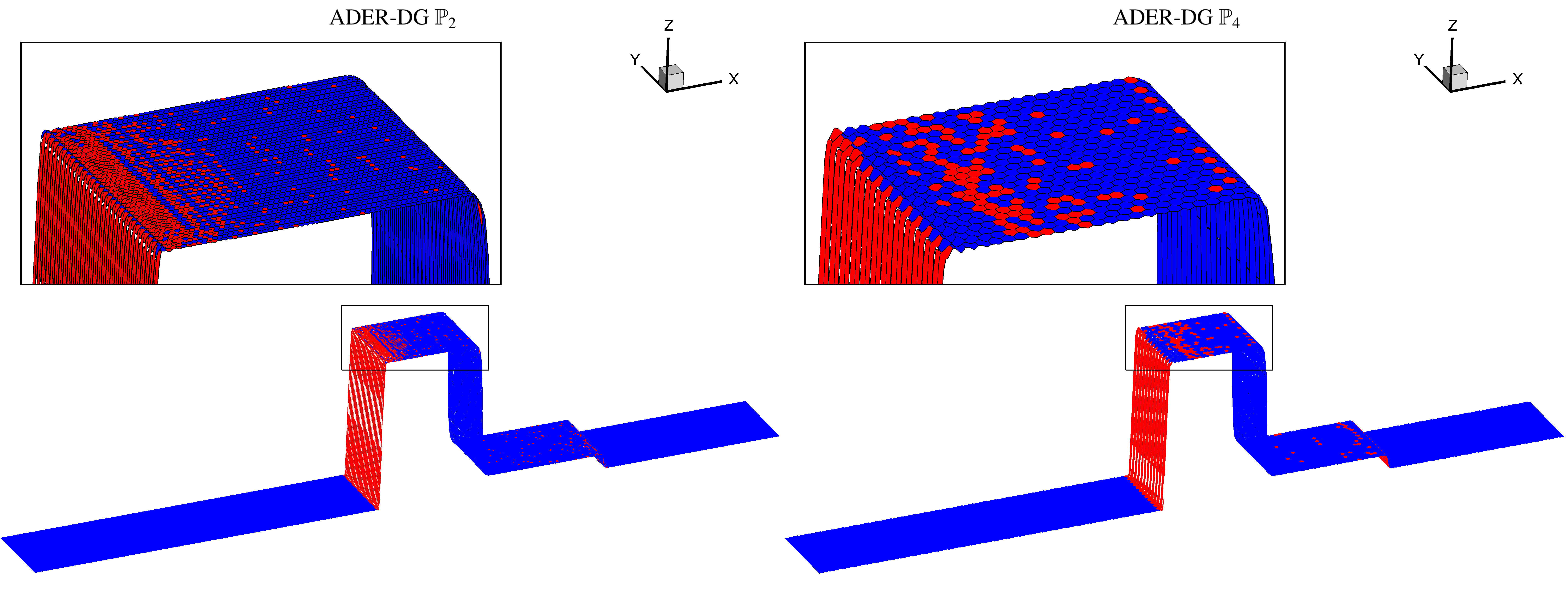}
		\caption{
		3D visualization of the \textit{a posteriori} subcell FV limiter's activation pattern for the solution of the \textit{2S problem} at final time $t = 0.4$, with the profile of the density $\rho$ represented along the $z$ axis. \textit{Troubled} elements, i.e. where the limiter is active, are shown in red, while unlimited elements are displayed in blue. Left panel: third order ADER-DG scheme. Right panel: fifth order ADER-DG scheme. The zoomed-in views highlight the underlying unstructured polygonal meshes.
		}
		\label{fig.RHD.Riemann3.2S.limiter}
	\end{center}	
\end{figure}

We run these tests over a rectangular computational domain $\Omega = [-0.5,0.5]\times[-0.05,0.05]$ with periodic boundary conditions in the $y$ direction, using the {\emph{harmonic gauge condition}} $H_{a} = 0$, and all damping parameters set to zero. 
We use the third order CWENO-FD scheme on $512 \times 8$ grid points with the WENO parameters of Eq.~\eqref{WENOwr} set to $r = 10$ and $\epsilon = 10^{-14}$, and the third order ADER-DG scheme on a mesh generated by $N_{x} = 500$ elements, combined with an \textit{a posteriori} second order ADER-ENO subcell finite volume limiter. The results are presented in Figures \ref{fig.RHD.Riemann1.2R}, \ref{fig.RHD.Riemann2.RS}, and \ref{fig.RHD.Riemann3.2S}, where a comparison with the exact profiles at the final time $t = 0.4$ shows that both classes of schemes effectively capture the challenging features of the proposed configurations.
To better illustrate the behavior of the ADER-DG scheme on unstructured meshes, in particular when combined with the presented limiting technique, we also report in Figure~\ref{fig.RHD.Riemann3.2S.limiter} a visualization of the \textit{a posteriori} subcell FV limiter's activation pattern for the solution of the \textit{2S problem} at the final time $t = 0.4$. We compare the results obtained with the third and fifth order versions of the scheme over $N_{x} = 500$ and $N_{x} = 200$ elements, respectively, in order to have the same number of degrees of freedom. We notice that, apart from some spurious triggering, the limiter is indeed used mainly in the vicinity of the strong shock waves, while the unlimited high-order DG scheme is retained in the smooth regions of the solution.

\subsection{Spherical Michel accretion}
\label{sec:spherical-accretion}
Accretion can be described as the accumulation of matter onto a compact object due to gravitational attraction, and it constitutes a well-studied problem in astrophysics.
The simplest classical model in accretion theory is the spherically symmetric flow of an infinite gas cloud onto a Newtonian point-mass, suggested by Hermann Bondi~\cite{bondi1952spherically}. A general-relativistic extension was proposed by F. Curtis Michel~\cite{michel1972accretion}, who considered the spherical accretion of a stationary flow onto a Schwarzschild spacetime (see~\cite{Rezzolla_book:2013} for a modern presentation of the problem).  
This represents one of the first studies of spherical accretion onto black holes, and since then, it has inspired several extensions accounting for additional physical effects.

In our work, the fluid is modeled as an ideal gas with adiabatic index $\gamma= 5/3$.
The key parameters are the critical radius, where the flow becomes supersonic, which we set to $r_c = 8$, and the critical density, defined as the density at the critical radius, that we choose as $\rho_c = 1/16$ according to~\cite{DelZanna2007}.
All the details on how to construct the initial condition are presented in~\cite{Anton06}.
It is essential to emphasize that the Michel solution assumes a static spacetime, and consequently rather than solving the full Einstein-Euler equations, we evolve only the Euler part while keeping the background metric fixed. This approach, i.e. the Cowling approximation, is justified when the mass accretion rate is small enough that the growth of the black hole mass, potentially originating from the accreted matter, can be safely neglected.

We perform the numerical tests in Kerr-Schild spheroidal coordinates~\cite{Komissarov04b} on a computational domain $(r,\theta) \in [1, 10]\times[\delta, \pi-\delta]$, with $\delta=0.05$.
The gauge source function $H_{a}$ is computed from the initial metric as $H_{a} = -\Gamma_{a}$, and we set all the damping parameters to zero. At the boundaries, we impose the Dirichlet initial data.

\begin{figure}[!htbp]
\begin{center}
\includegraphics[width=0.8\textwidth, keepaspectratio]{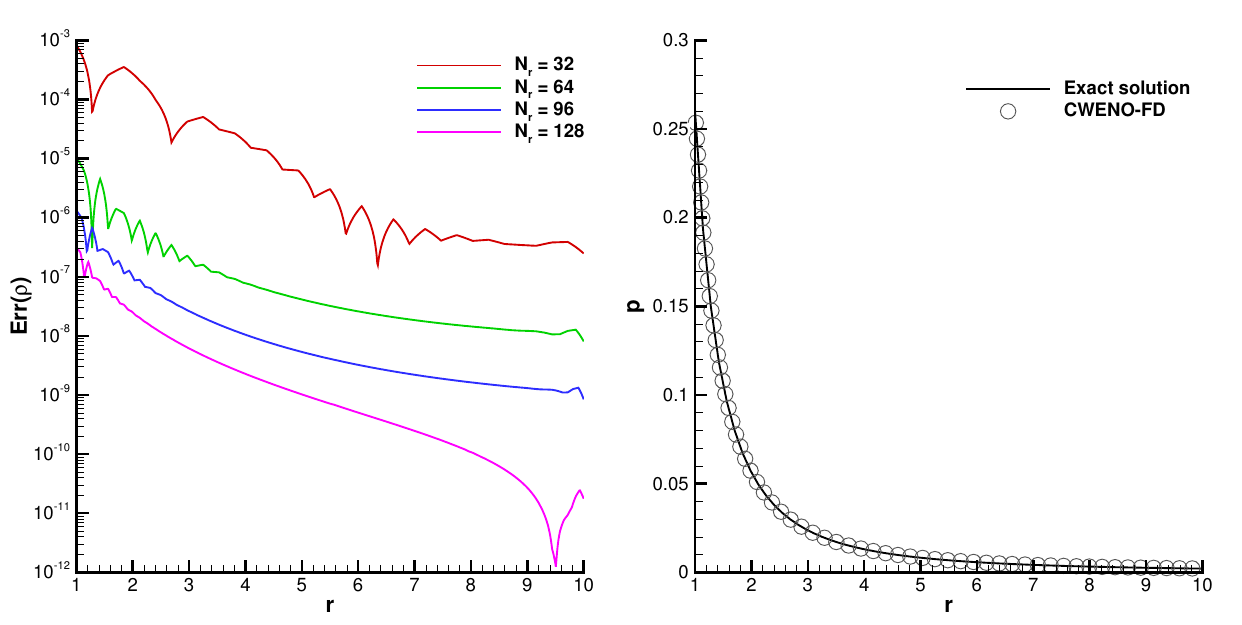}
\caption{
Spherical accretion onto Schwarzschild black hole in 2D spheroidal Kerr-Schild coordinates using the fifth order CWENO-FD  scheme. Left panel: absolute errors (computed with respect to the exact solution) for the density $\rho$ at final time $t = 1000$ while increasing the number $N_{r}$ of radial grid points. Right panel: comparison of the exact and numerical solution for the pressure $p$ at final time $t = 1000$ on a one dimensional cut at $\theta = \pi/2$.
}
\label{fig.SphericalAccretion.CWENO4.64}
\end{center}	
\end{figure}

We show the results obtained using the fifth order CWENO-FD scheme in Figure~\ref{fig.SphericalAccretion.CWENO4.64}, where in the left panel we demonstrate mesh convergence by plotting the absolute errors for the density $\rho$ at the final time $t=1000$ for different radial resolutions. Additionally, in Tables~\ref{tab:SphericalAccretion.FDToy.convergence} and \ref{tab:SphericalAccretion.AleVoronoi.convergence} we present a convergence analysis at time $t = 1$ for the CWENO-FD and ADER-DG schemes, respectively.

\begin{table}[!htbp]
\centering
\begin{tabular}{lcccc}
\toprule
\textbf{Scheme} & $\bm{N_{r}\times N_{\theta}}$ & $\bm{L^2}$ \textbf{Error} $\bm{\tilde{D}}$ & $\bm{\mathcal{O}(\tilde{D})}$ & \textbf{Theoretical}\\
\midrule
CWENO-FD $\mathbb{P}_{4}$ & $128\times64$  & 3.1814E-08 & -\\
 & $256\times128$ & 1.4794E-09 & 4.43 & 5.00\\
 & $512\times256$ & 5.7317E-11 & 4.69  \\
\midrule
CWENO-FD $\mathbb{P}_{6}$ & $128\times64$ & 7.8714E-10 & -  \\
 & $192\times96$ & 7.2139E-11 & 5.89 & 7.00 \\
 & $256\times128$ & 1.2307E-11 & 6.15 \\
\bottomrule
\end{tabular}
\caption{Numerical convergence results for the spherical accretion onto Schwarzschild black hole in 2D spheroidal Kerr-Schild coordinates at final time $t=1$ for the conserved variable $\tilde{D} \coloneqq  \sqrt{\gamma}D = \sqrt{\gamma}\rho W$ using the CWENO-FD schemes.
}
\label{tab:SphericalAccretion.FDToy.convergence}
\end{table}
\begin{table}[!htbp]
\centering
\begin{tabular}{lccccc}
\toprule
\textbf{Scheme} & $\bm{N_{x}}$ & $\bm{h}$ & $\bm{L^2}$ \textbf{Error} $\bm{\tilde{D}}$ & $\bm{\mathcal{O}(\tilde{D})}$ & \textbf{Theoretical}\\
\midrule
ADER-DG $\mathbb{P}_{1}$ & $80$ & 4.2075E-02 & 6.1397E-03  & - \\
 & $120$ & 2.7900E-02 & 2.6688E-03 & 2.03 & 2.00\\
 & $160$ & 2.1037E-02 & 1.5310E-03 & 1.97 \\
\midrule
ADER-DG $\mathbb{P}_{2}$ & $40$ & 8.2836E-02  & 5.9163E-04 & -  \\
 & $80$ & 4.2075E-02  & 8.2674E-05 & 2.91 & 3.00\\
 & $120$ & 2.7900E-02 & 2.3864E-05 & 3.02 \\
\midrule
ADER-DG $\mathbb{P}_{3}$ & $40$ & 8.2836E-02 & 1.4810E-05  & - \\
 & $80$ & 4.2075E-02 & 9.9158E-07 & 3.99 & 4.00\\
 & $120$ & 2.7900E-02 & 1.8447E-07  & 4.09 \\
 \midrule
ADER-DG $\mathbb{P}_{4}$ & $20$ & 1.6567E-01 & 4.5139e-06 & - \\
 & $40$ & 8.2836E-02& 1.9273E-07  & 4.55 & 5.00\\
 & $80$ & 4.2075E-02 & 8.9232E-09  & 4.54 \\
\bottomrule
\end{tabular}
\caption{Numerical convergence results for the spherical accretion onto Schwarzschild black hole in 2D spheroidal Kerr-Schild coordinates at final time $t=1$ for the conserved variable $\tilde{D} \coloneqq  \sqrt{\gamma}D = \sqrt{\gamma}\rho W$ using the ADER-DG schemes.
}
\label{tab:SphericalAccretion.AleVoronoi.convergence}
\end{table}

\subsection{Non-rotating neutron star in equilibrium}
\label{sec:tov}
Neutron stars are extremely dense astrophysical objects that result from the gravitational collapse of massive star cores.  Even though formulating comprehensive models of their structure is quite complex, there exist some relatively simple standard setups that are commonly used to test the ability of numerical relativity codes to evolve the full Einstein-Euler system, while also giving a first insight into real astrophysical applications.
A representative one is the evolution of an {\emph{equilibrium non-rorating neutron star} over long timescales.
We therefore consider the Tolman-Oppenheimer-Volkoff (TOV) equations~\cite{Oppenheimer39b, Tolman}
\begin{align}
\label{eq:tov1}
\frac{d m}{d r} &= 4\pi r^2 (\rho h - p) \,,\\
\label{eq:tov2}
\frac{dp}{dr}   &= - \frac{\rho h(m+4\pi r^3 p)}{r(r-2m)}\,,\\
\label{eq:tov3}
\frac{d \phi}{d r} &= -\frac{1}{\rho h}\frac{dp}{dr} \,,
\end{align}
where $m(r)$ is the mass enclosed within the radius $r$, $\phi$ is the unknown metric function in the line element
\begin{equation}\label{eq.symm.metric}
ds^2 = -e^{2\phi} dt^2 + e^{2\psi} dr^2 + r^2 d\Omega^2\,,
\end{equation}
while $e^{-2\psi}=1-\frac{2m}{r}$.
When coupled with an equation of state (EoS), the TOV system completely determines the physical properties of a spherically symmetric body in equilibrium, such as the internal pressure and density profiles, as well as its total mass and radius.
Here we consider the polytropic EoS $p= K\,\rho^\gamma$, and we choose the model parameters according to~\cite{Font2002}, meaning a central rest mass density $\rho_c=1.28\times 10^{-3}$, $K=100$, and $\gamma=2$. 
We then integrate numerically the ODEs~\eqref{eq:tov1}--\eqref{eq:tov3} with a tenth order accurate discontinuous Galerkin scheme~\cite{ADERNSE}, obtaining a total mass $M=1.4\,M_{\odot}$ and a radius $R=9.585 \,M_{\odot}= 14.15\,km$.
As described in~\cite{bugner}, the standard approach is to introduce the isotropic radius $\bar{r} = \bar{r}(r)$, and to apply the change of coordinates
\begin{equation}
\frac{d\bar{r}}{\bar{r}}=\left(1-\frac{2m}{r}\right)^{-1/2}\frac{dr}{r}, 
\end{equation}
which makes the spatial part of the metric~\eqref{eq.symm.metric} conformally flat, i.e.
\begin{equation}
ds^2 = -e^{2\phi} dt^2 + e^{2\bar{\psi}}( d\bar{r}^2 + \bar{r}^2 d\Omega^2) =-e^{2\phi} dt^2 + e^{2\bar{\psi}}(d\bar{x}^2+d\bar{y}^2+d\bar{z}^2  )\,,
\end{equation}
where $e^{2\bar{\psi}} {\bar{r}}^{2} = r^{2}$.
In the exterior of the star, the spacetime is that of a Schwarzschild solution produced by a mass $M$, as a consequence of Birkoff’s theorem, while all the hydrodynamic variables vanish. We underline that, contrary to what is usually done, we do not impose any artificial atmosphere since we use the filter for the conversion from the conserved to the primitive variables introduced in~\cite{DumbserZanottiGaburroPeshkov2023}.
To conclude the description of the construction of the solution, we recall that since we need the initial data at arbitrary $\bar{r}$, interpolation is needed.

For our simulations, we consider a three-dimensional domain $\Omega = [-15, 15]^{3}$ discretized with $128^{3}$ points, and we use the fifth order version of the CWENO-FD schemes. The gauge source function $H_{a}$ is computed from the initial data through $H_{a} = -\Gamma_{a}$, and we set the damping parameters to $\gamma_{0} = 0.1, \gamma_{1} = -1, \gamma_{2} = 0.5$.
We add a small perturbation of amplitude $10^{-4}$ to the initial pressure $p$, so that $\tilde{p} = p\,(1 + 10^{-4})$, and we compare the results obtained with and without well-balancing.  More precisely, in the left panel of Figure~\ref{fig.TOV.rhomax.WB.vs.noWB} we show the time evolution of the normalized central density up to final time $t = 1000$, from which it is evident that the well-balanced approach leads to a much more stable dynamics, while in the right panel of the same figure we further extend the simulation up to $t = 2000$ for the well-balanced case only, in order to highlight its long term robustness. Finally, in Figure~\ref{fig.TOV.constraints.WB.vs.noWB}, we show the time evolution of the constraints.
\begin{figure}[!htbp]
\begin{center}
\includegraphics[width=0.8\textwidth, keepaspectratio]{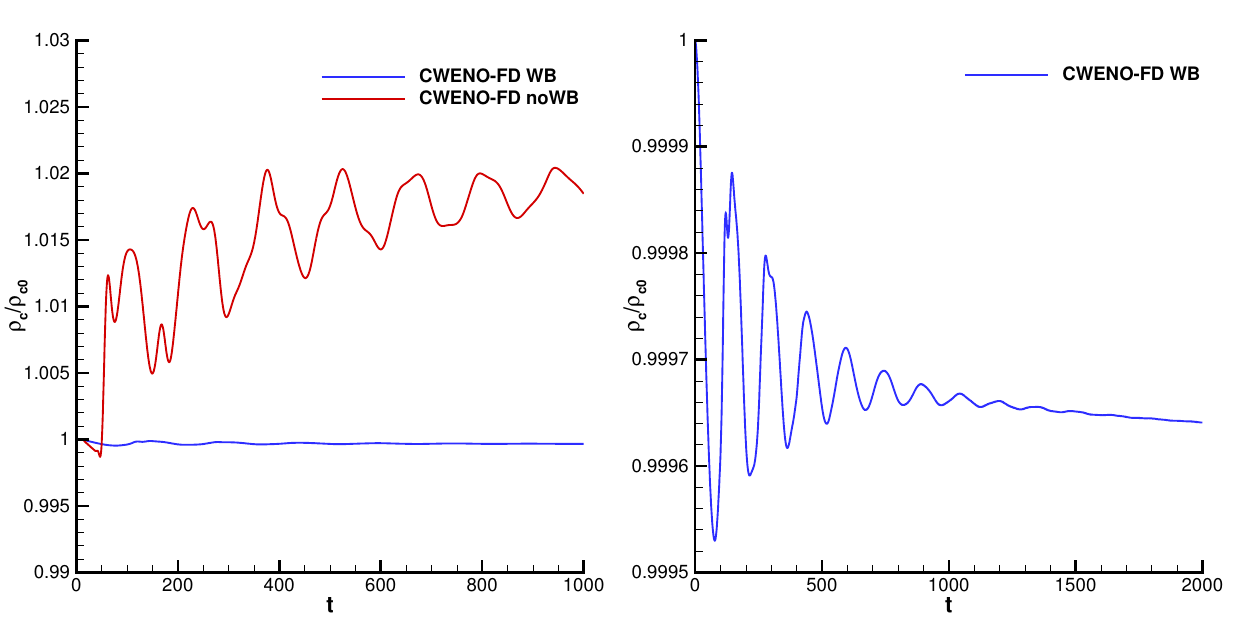}
\caption{Time evolution of the normalized central density for the perturbed non-rotating neutron star solved with the fifth order CWENO-FD scheme. Left panel: comparison of the well-balanced (WB) and non well-balanced (noWB) versions up to the final time $t = 1000$. Right panel: profile up to the final time $t = 2000$ only in the well-balanced case. 
}
\label{fig.TOV.rhomax.WB.vs.noWB}
\end{center}	
\end{figure}  
\begin{figure}[!htbp]
\begin{center}
\includegraphics[width=0.45\textwidth, keepaspectratio]{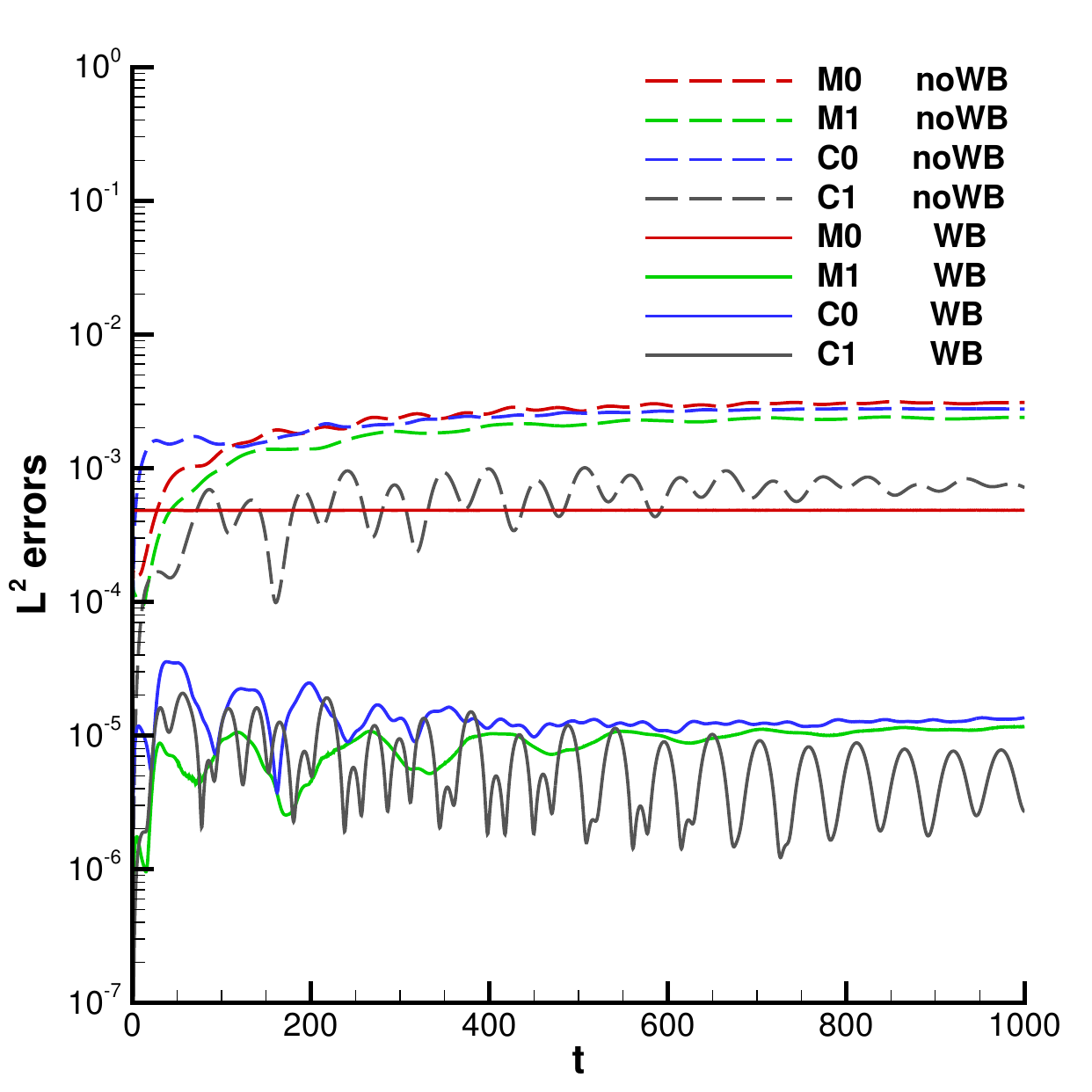}
\caption{Time evolution of the constraints for the perturbed non-rotating neutron star. Comparison between the well-balanced (WB) and non well-balanced (noWB) fifth order CWENO-FD schemes.
}
\label{fig.TOV.constraints.WB.vs.noWB}
\end{center}	
\end{figure}

\section{Conclusions}
\label{sec.conclusions}
In this paper, we have tested two distinct families of high order numerical schemes for nonlinear hyperbolic systems of balance laws applied to the Generalized Harmonic (GH) formulation of the Einstein equations first introduced in~\cite{Lindblom2006}: a CWENO-FD scheme on Cartesian grids in two and three space dimensions, and an ADER-DG scheme on unstructured polygonal meshes in two space dimensions. Both approaches were endowed with a well-balancing feature, enabling the exact discrete preservation of known stationary solutions. The successful reproduction of a wide variety of standard vacuum benchmarks confirmed the robustness and accuracy of our schemes when applied to numerical relativity problems.
Moreover, the inclusion of matter evolution through the general relativistic Euler equations allowed us to further validate our methods also in a more complete regime, namely coupling the spacetime and matter parts, through tests such as the spherical accretion onto a Schwarzschild black hole and the non-rotating neutron star in equilibrium. 
In addition, relativistic Riemann problems were considered in order to assess the ability of our schemes to correctly capture strong shock waves and discontinuities.
For the test cases involving the proposed well-balancing techniques, a rather natural choice was to focus on small perturbations around stationary equilibria.
However, for more relevant physical applications of these tools we mention the stability analysis of black hole horizons \cite{Aretakis2015,Agrawal2026}, as well as oscillation modes of neutron stars \cite{servignat2025roxas,Jaraba2026,Zhao2022}.
More dynamic simulations, such as binary neutron star mergers, as well as further developments concerning mycrophysics and radiation transport will be the subject of future research. 
Overall, the numerical results obtained in this paper establish a solid basis for the next stages of development of this work, for which it is appropriate to offer a few additional comments. A first natural extension consists in generalizing the ADER-DG scheme to a complete 3D framework on unstructured meshes, ultimately supporting the study of more realistic astrophysical sources beyond axisymmetric systems. Furthermore, the integration of moving meshes within an Arbitrary-Lagrangian-Eulerian (ALE) framework with topology changes~\cite{Gaburro2020Arepo} represents another key step towards accurate simulations of dynamical spacetimes, such as those encountered in rotating binary systems. We recall that, within the GH formulation, neutron star and black hole mergers have recently been simulated successfully with high order DG schemes in~\cite{Deppe2024} and~\cite{lovelace2025simulating}, respectively. The first of these two scenarios, in fact, constitutes the main target we aim to address in our future works, potentially including the extraction of gravitational waves.

\section{Acknowledgments}

All the authors of this paper are members of the Gruppo Nazionale Calcolo Scientifico-Istituto Nazionale di Alta Matematica (GNCS-INdAM) and would like to acknowledge the CINECA award under the ISCRA initiative, for the availability of high-performance computing resources and support.

M.~Dumbser and O.~Zanotti acknowledge the support received via the  Departments of Excellence Initiative 2018--2027 attributed to DICAM of the University of Trento (grant L. 232/2016).
M.~Dumbser also acknowledges the funding received via the Fondazione Caritro under the project SOPHOS. 

E.~Gaburro and S.~Muzzolon gratefully acknowledge the support received from the European Union 
with the ERC Starting Grant \textit{ALcHyMiA} (grant agreement No. 101114995). 
Views and opinions expressed are however those of the authors only and
do not necessarily reflect those of the European Union or the European Research
Council. Neither the European Union nor the granting authority can be held
responsible for them.



\bibliographystyle{plain}
\bibliography{./referencesZ4.bib}

\end{document}